\documentclass[11pt]{article}

\usepackage{amsmath,amsfonts,amssymb,amsthm,mathtools,bbm, geometry, cancel, tikz}

\usetikzlibrary{matrix}

\usepackage{graphicx, multicol}

\DeclareMathAlphabet{\mathboondoxfrak}{U}{BOONDOX-frak}{m}{n}

\usepackage{amssymb}

\newtheorem{thm}{Theorem}[section]

\newtheorem{dfn}[thm]{Definition}

\newtheorem{lem}[thm]{Lemma}

\newtheorem{prop}[thm]{Proposition}

\newtheorem{remark}[thm]{Remark}

\newtheorem{cor}[thm]{Corollary}

\newtheorem{ex}[thm]{Example}

\newtheorem{question}[thm]{Question}

\def\bq{\begin{question}}
\def\bt{\begin{thm}}
\def\bp{\begin{prop}}
\def\blem{\begin{lem}}
\def\bd{\begin{dfn}}
\def\br{\begin{remark}}
\def\bc{\begin{cor}}
\def\bex{\begin{ex}}
\def\beqs{\begin{eqnarray*}}
\def\beq{\begin{eqnarray}}
\def\bi{\begin{itemize}}

\def\eq{\end{question}}
\def\et{\end{thm}}
\def\ep{\end{prop}}
\def\elem{\end{lem}}
\def\ed{\end{dfn}}
\def\er{\end{remark}}
\def\ec{\end{cor}}
\def\eex{\end{ex}}
\def\eeqs{\end{eqnarray*}}
\def\eeq{\end{eqnarray}}
\def\ei{\end{itemize}}

\def\no{\noindent}
\def\ds{\displaystyle}

\def\ov{\overline}
\def\r{\rangle}
\def\l{\langle}

\def\H{{\cal H}}
\def\K{{\cal K}}
\def\B{{\cal{B}}}

\def\Hp{{\cal H}_\pi}

\def\w*{w^*-w^*}

\def\F{{\cal F}}

\def\ra{\rightarrow}

\def\vp{\varphi}
 
\def\bs{\backslash}

\def\E{{\cal E}}

\def\K{{\mathcal K}}

\makeatletter
\def\moverlay{\mathpalette\mov@rlay}
\def\mov@rlay#1#2{\leavevmode\vtop{%
   \baselineskip\z@skip \lineskiplimit-\maxdimen
   \ialign{\hfil$\m@th#1##$\hfil\cr#2\crcr}}}
\newcommand{\charfusion}[3][\mathord]{
    #1{\ifx#1\mathop\vphantom{#2}\fi
        \mathpalette\mov@rlay{#2\cr#3}
      }
    \ifx#1\mathop\expandafter\displaylimits\fi}
\makeatother

\newcommand{\cupdot}{\charfusion[\mathbin]{\cup}{\cdot}}
\newcommand{\bigcupdot}{\charfusion[\mathop]{\bigcup}{\cdot}}

\usepackage[all]{xy}

\usetikzlibrary{arrows.meta,positioning}
\usetikzlibrary{cd}

\definecolor{c1}{HTML}{1F77B4} 
\definecolor{c2}{HTML}{D62728} 
\definecolor{c3}{HTML}{2CA02C} 
\definecolor{c4}{HTML}{FF7F0E} 
\definecolor{c5}{HTML}{9467BD} 
\definecolor{c6}{HTML}{8C564B}

\makeatletter
\providecommand*{\cupdot}{%
  \mathbin{%
    \mathpalette\@cupdot{}%
  }%
}
\newcommand*{\@cupdot}[2]{%
  \ooalign{%
    $\m@th#1\cup$\cr
    \sbox0{$#1\cup$}%
    \dimen@=\ht0 %
    \sbox0{$\m@th#1\cdot$}%
    \advance\dimen@ by -\ht0 %
    \dimen@=.5\dimen@
    \hidewidth\raise\dimen@\box0\hidewidth
  }%
}

\providecommand*{\bigcupdot}{%
  \mathop{%
    \vphantom{\bigcup}%
    \mathpalette\@bigcupdot{}%
  }%
}
\newcommand*{\@bigcupdot}[2]{%
  \ooalign{%
    $\m@th#1\bigcup$\cr
    \sbox0{$#1\bigcup$}%
    \dimen@=\ht0 %
    \advance\dimen@ by -\dp0 %
    \sbox0{\scalebox{2}{$\m@th#1\cdot$}}%
    \advance\dimen@ by -\ht0 %
    \dimen@=.5\dimen@
    \hidewidth\raise\dimen@\box0\hidewidth
  }%
}
\makeatother

\def\S{{\cal S}}

\def\R{\mathbb{R}}

\def\fza{F_z^\alpha} 
\def\kza{K_z^\alpha} 
\def\bza{\beta_z^\alpha} 
\def\fwa{F_w^\alpha} 
\def\kwa{K_w^\alpha} 
\def\bwa{\beta_w^\alpha} 

\def\preq{\preceq_Q}
\def\TG{{\cal T}(G)}
\def\TnqG{{\cal T}_{\text{nq}}(G)}
\def\TnqA{{\cal T}_{\text{nq}}^\A}
 
\def\ABSG{A_\B^*(G)} 
\def\taunq{{\tau_{\text{nq}}}}
\def\signq{{\sigma_{\text{nq}}}}
\def\D{{\cal D}} 
\def\ADSG{A_\D^*(G)}
\def\nq{{\text{nq}}} 
\def\HD{\text{HD}}
\def\PHD{\text{PHD}}

\def\N{\mathbb{N}}
\def\A{{\mathfrak A}} 
\def\E{{\cal E}}
\def\F{{\cal F}}

\begin{document}

\title{(Generalized) Spine Subalgebras  of Fourier--Stieltjes algebras and their Homomorphisms}
\author{Nico Spronk$^*$,  Ross Stokke\footnote{These authors received partial support  from NSERC grants.}
   \ and Aasaimani Thamizhazhagan\footnote{
This author received financial support from the Department of Science and Technology, India, under the FIST program (Ref. SR/FST/MS-I/2024/173).}  }
\date{}
\maketitle

\begin{abstract}{\small    For any upper semilattice $\D$ of locally precompact topologies on a locally compact group $G$,  we define an associated generalized spine subalgebra  $A^*_\D(G)$ of the Fourier--Stieltjes algebra $B(G)$. We show that  $A^*_\D(G)$ is a semilattice-graded $\ell^1$-direct sum of maximal copies of Fourier algebras and we identify its spectrum as a semilattice of groups. We build a collection of examples of generalized spine algebras over whose spectra we  exhibit fine control.  We define notions of compatible fusions of homomorphisms and affine maps, and use these definitions to  characterize all completely positive, completely contractive and, when $G$ is amenable, all completely bounded homomorphisms from a generalized spine algebra $A^*_\D(G)$ to a Fourier--Stieltjes algebra $B(H)$. These results are new, even when $A^*_\D(G)$ is the full spine algebra $A^*(G)$ and even when $G$ and $H$ are abelian. We provide examples illustrating the scope of our theorems. 

\smallskip

 }
\end{abstract} 

\section{Introduction}

\subsection{Background and plan} 
For a locally compact group, $G$, the primary objects in abelian and noncommutative harmonic analysis are its group and measure algebras $L^1(G)$ and $M(G)$ with convolution product, and, dually, its Fourier and Fourier--Stieltjes algebras $A(G)$ and $B(G)$ with pointwise-defined product.  When $G$ is abelian, $A(G) \cong L^1(\widehat{G})$ and $B(G) \cong M(\widehat{G})$, through the Fourier--Stieltjes transform on $\widehat{G}$. Integral to our understanding of any commutative Banach algebra $A$ is its Gelfand spectrum $\Delta(A)$, the space of nonzero multiplicative linear functionals on $A$. With the Gelfand spectra of $A(G)$ and $B(G)$, we typically see a dichotomy of complexity:  $\Delta(A(G))$, which can be identified with $G$, is fully accessible;  meanwhile, though it can sometimes be described, in the words of E. Kaniuth and A.T.-M. Lau, ``when $G$ is a noncompact locally compact abelian group, according to common understanding, the spectrum of $B(G) = M(\widehat{G})$ is an intractable object" \cite[Section 2.9]{Kan-Lau}.   This dichotomy is  one of the   reasons why $A(G)$ is more accessible than $B(G)$, though we note that Walter \cite{Wal2} and others have provided significant insights into the structure of $\Delta(B(G))$.  

  When  $G$ is abelian, the spine algebra $L^*(G)$ of $M(G)$ was independently defined and studied by J. Inoue \cite{Inoue} and J.L. Taylor, who introduced the terminology \cite{Tay}.  M. Ilie and the first-named author later defined the spine $A^*(G)$ of $B(G)$ for any locally compact group $G$ \cite{Il-Spr1}. Through the Fourier-Stieltjes transform, $A^*(G) = L^*(\widehat{G})$ when $G$ is abelian. 
We let $\TG$ denote the join semilattice  of all locally precompact topologies on $G$; so,   $\tau \in \TG$ if it is  the weak topology  induced by  a continuous dense-range homomorphism $\eta_\tau$ of $G$ into a locally compact group  $G_\tau$. For $\tau$ in $\TG$, $A_\tau(G): = A(G_\tau) \circ \eta_\tau$ is an isometric copy of the Fourier algebra $A(G_\tau)$ in $B(G)$, and the spine of $B(G)$ is defined in \cite{Il-Spr1} as $A^*(G):= \ov{\sum_{\tau \in \TG} A_\tau(G)}$,  the closed linear span in $B(G)$ of all the subalgebras $A_\tau(G)$.  This canonical object has spectrum that is usually of intermediate complexity between the spectra of $A(G)$ and $B(G)$: it is a semigroup  and a  union of groups graded over a  semilattice determined by $\TG$, i.e., $\Delta(A^*(G))$ is a semilattice of groups \cite[Theorem 4.1]{Il-Spr1}. Though $\Delta(A^*(G))$ can be complicated, it is not  inaccessible.  A major achievement in Taylor's work on abelian groups was  showing that the invertible elements $M(G)$ must lie in its more tractable spine algebra $L^*(G)$, a result extended to $B(G)$ and $A^*(G)$  for some classes of non-abelian groups by the third-named author \cite{Tha}.

The problem of describing homomorphisms $\vp: A  \ra B(H)$ demonstrates the importance of understanding $\Delta(A)$ when $A$ is a closed translation-invariant subalgebra of $B(G)$.  In the case that $A$ is the Fourier algebra $A(G)$, this is part of an old problem that was solved for abelian groups by Paul Cohen in 1960, but remains open in general. For non-abelian groups and $A=A(G)$, having characterized every such (completely) positive, (completely) contractive and, for $G$ amenable, completely bounded homomorphism,   M. Ilie \cite{Il}, Ilie--Spronk \cite{Il-Sp}, M. Daws \cite{Daws}  and H.L. Pham \cite{Pham} have, among others, made significant contributions to this problem. (The dual problem, which asks for descriptions of homomorphisms between the convolution algebras $L^1(G)$ and $M(H)$ for nonabelian groups, is considered for contractive and positive homomorphisms in \cite{Gre, KSSY, Sto11}, for example.) 

It seems equally natural to study homomorphisms $\vp: A \ra B(H)$ when $A$ is another closed translation-invariant subalgebra of $B(G)$ \cite{Il-St, Pham, Sto21, Wal1}, however, since any such homomorphism is uniquely determined  by an open subset $E$ of $H$ and a continuous map $\alpha: E \ra \Delta(A)$ --- we write $\vp = j_\alpha$ --- an understanding of $\Delta(A)$ (or at least the image of $\alpha$ in $\Delta(A)$) is needed to tackle this problem.  Hence, it is particularly  natural to study the homomorphisms $\vp: A^*(G) \ra B(H)$ and, when $G$ and $H$ are abelian, Inoue characterized every such homomorphism in \cite{Inoue}. (As observed in \cite{Sto21}, e.g., see Examples 2.14 and 4.16 therein, there is an issue with the descriptions in \cite{Il-Spr1} of these homomorphisms:  $\alpha$ may fail  to map into a single subgroup of $\Delta(A^*(G))$, even when $\vp=j_\alpha$ is completely positive.) 
However, while the spectrum of the full spine algebra $A^*(G)$ is a semilattice of groups that was  explicitly computed for several groups in \cite{Il-Spr1}, these instances do not alone provide the breadth of examples we needed to answer some of our questions about homomorphisms $\vp: A \ra B(H)$. 

Thus motivated, taking $\D$ to be any join subsemilattice of $\TG$,   in Section 3 
we define the associated generalized spine algebra (GSA) to be the subalgebra $\A = A_\D^*(G) = \ov{ \sum_{\tau \in \D} A_\tau(G) }$ of $A^*(G)$. We show --- often by employing results from \cite{Il-Spr1} -- that most of the theory developed in \cite{Il-Spr1}  for $A^*(G)$ carries over to generalized spine algebras $\A$. For example, there exist maximal subalgebras $A_\tau(G)$ for $\tau$ in an associated  semilattice  $(\B, \vee_\A)$  of locally precompact topologies   such that  $\A = \ell^1-\oplus_{\tau \in \B}  A_\tau^*(G)$ with $A_\tau(G) A_\sigma(G) \subseteq A_{\tau \vee_\A \sigma}(G)$,  and we show that $\Delta(\A)$ is a semilattice of groups. In Section 4,  for certain groups $G$ and any upper subsemilattice $\E$ of  $({\cal P}(\N), \cup)$, we build  a GSA, $A_\E^*(G) = \ell^1-\oplus_{A \in \E} A_{\tau_A}(G)$, graded over $(\E, \cup)$.  When studying homomorphisms $\vp: \A \ra B(H)$,  the most interesting aspect of the semilattice of groups $\Delta(\A)$  is the indexing semilattice  itself and, by selectively choosing $\E$, we obtain GSAs whose spectra $\Delta(A^*_\E(G))$  have indexing semilattices ${\cal L}$ over which  we can impose considerable control: for example, we  can ensure that ${\cal L}$ is any finite Boolean subalgebra, an infinite chain, or  the full lower semilattice $({\cal P}(\N), \cap)$.     We anticipate that the constructions in Section 4 will therefore be interesting to researchers interested in abelian or noncommutative harmonic analysis, or the general theory of Gelfand spectra.

In Section 5, we extend the definitions from \cite{Sto21} of compatible fusions of homomorphisms to the context needed herein, and explore some properties of these maps.  Letting $\A$ be any GSA,  in Section 6 we characterize all completely positive,   completely contractive and, for $G$ amenable,  completely bounded  homomorphisms $\vp: \A \ra B(H)$.  These results are new,  even when $\A$ is the full spine algebra $A^*(G)$ and even when  $G$ and $H$ are abelian, though for abelian groups our description of the (completely) bounded homomorphisms is related to Inoue's.  For example, Theorems \ref{MainCPCCHomomThm} and \ref{CPCCHomomThmsNoSemigp} show that  $\vp: \A \ra B(H)$ is a completely positive (contractive) homomorphism if and only if $\vp = j_\alpha$ where $\alpha$ is a  continuous fusion of compatible homomorphisms (affine maps) mapping an open subgroup (coset) of $H$ into $\Delta(\A)$. 
Answering some natural questions arising from \cite{Sto21}, (see Remark \ref{Scope of alpha for CP Homs Remark}.2),  
 Example \ref{MainCPHomExample} shows that such a fusion of homomorphisms can map nontrivially into every subgroup of $\Delta(\A)$, where $\Delta(\A)$ is faithfully graded over $({\cal P}(\N), \cap)$, an uncountable lower semilattice containing infinite chains and anti-chains.

\subsection{Basic definitions} 

Throughout this paper, $G$ and $H$ are  locally compact groups.   Unless stated otherwise, we follow the notation in \cite{Ars, Eym, Kan-Lau}, where the reader can find the  properties of the Fourier algebra $A(G)$ and the Fourier--Stieltjes algebra $B(G)$ used herein. If   $\pi$  is  a continuous unitary representation of $G$ on a Hilbert space, $\H_\pi$, 
$$ \phi^\pi_{\xi, \eta}(s):=   \l \pi (s) \xi | \eta \r \qquad  (s \in
G, \ \xi, \eta \in \H)$$ 
are its coefficient functions and the Fourier space $A_\pi$  is the closed linear span in $B(G)$  of all coefficient functions of $\pi$.   We use ${\cal L}(\H_\pi)$ to denote the space of bounded linear operators on $\H_\pi$ and $VN_\pi$ is the  von Neumann subalgebra of ${\cal L}(\Hp)$ generated by $\pi(G)$; $VN_\pi$  can be identified with the dual of $A_\pi$ through the pairing 
$\l T, \phi^\pi_{\xi, \eta} \r = \l T \xi | \eta\r$.  The Fourier spaces $A_\pi$ are precisely the closed translation-invariant subspaces of $B(G)$.   Letting $\{ \lambda_G, L^2(G)\}$  and  $\{ \omega_G, \H_{\omega_G}\}$  be  the left regular and universal representations of $G$, respectively, $A(G) = A_{\lambda_G}$, $B(G) = A_{\omega_G}$, and we write    $VN(G)$  and $W^*(G)$   for $VN_{\lambda_G}$ and  $VN_{\omega_G}$, respectively.    

Since $A_\pi$ is the predual of a von Neumann algebra, it has a canonical operator space structure that  agrees with its subspace operator space structure inherited from $B(G) = W^*(G)_*$.   When it is a subalgebra of $B(G)$,  $A_\pi$ is a completely contractive Banach algebra.  
A map $\vp: A _\pi \ra B(H)$ is (completely) positive when its dual map $\vp^*: W^*(H) \ra VN_\pi$ is (completely) positive.  Standard references on the theory of operator spaces and completely bounded maps are  \cite{Eff-Rua} and   \cite{Pau}.

 A subset $E$ of $H$ is a coset of some subgroup of $H$ exactly when $E E^{-1}E = E$, and a map $\alpha: E \rightarrow G$ is   affine  if for any
$x,y,z \in E$, $\alpha(xy^{-1}z) =
\alpha(x)\alpha(y)^{-1}\alpha(z)$ \cite{Il}.  We let $\Omega(H)$ denote the ring of sets generated by the open cosets of $H$ and call a map $\alpha : Y \rightarrow G$  
piecewise affine  if
 there are pairwise disjoint sets $Y_1,...,Y_n \in \Omega(H)$ such that $Y=\bigcupdot_{i=1}^{n}{Y_i}$ and each  $\alpha\big{|}_{Y_i}$ has an affine
  extension, $\alpha_i$, mapping a coset containing $Y_i$  into $G$ \cite{Coh, Il}. \rm

  We revisit some ideas from Section 2 of \cite{Il-Spr1}; details are  found therein or readily verified.  Let $\tau_G$   denote the original locally compact group topology on $G$. A topology $\tau$ on $G$ is locally precompact if there is a locally compact group $(G_\tau, \tau_{G_\tau})$ and a continuous homomorphism $\eta_\tau: G \ra G_\tau$ such that $\tau$ is the weak topology on $G$ induced by $\eta_\tau$, i.e., if $$\tau = \eta_\tau^{-1}(\tau_{G_\tau}):=  \{ \eta_\tau^{-1}(V): V \text{ is open in } G_\tau\}.$$ Assuming without loss of generality that $\eta_\tau$ has dense range, the pair $(G_\tau, \eta_\tau)$ is unique up to an intertwining topological group isomorphism and $(G_\tau, \eta_\tau)$ is called a representation of $\tau$.  
We let $\TG$   denote  the set of locally precompact topologies on $G$.  (In \cite{Spr}, where non-locally compact topologies are also studied, ${\cal T}_{\text{lc}}(G)$ is used.)  For $\tau $ in $\TG$, observe that net convergence in $\tau$ is described by  \beq \label{NetCvgceTauEqn} t_\alpha \ra t \text{ in } \tau \quad  \text{ if and only if } \quad \eta_\tau(t_\alpha) \ra  \eta_\tau(t) \text{ in } G_\tau;\eeq  so, $\tau$ is a topological group topology on $G$.  Also, note that $\tau$ is Hausdorff exactly when $\eta_\tau$ is an injection.  If  $\sigma, \tau \in \TG$, then their least upper bound taken in the lattice of topologies of $G$, $\sigma \vee \tau$,  belongs to $\TG$ and has representation $(G_{\sigma \vee \tau}, \eta_{\sigma \vee \tau})$ where $G_{\sigma \vee \tau}$ is the closure of $\{(\eta_\sigma(s), \eta_\tau(s)): s \in G\}$ in $G \times G$ and $\eta_{\sigma\vee\tau}(s) = (\eta_\sigma(s), \eta_\tau(s))$. Thus, $(\TG, \vee)$ is an upper semilattice.  

When $\sigma \subseteq \tau$ in ${\cal T}(G)$, there is a unique continuous  homomorphism $\eta^\tau_\sigma$ such that $\eta^\tau_\sigma  \circ \eta_\tau = \eta_\sigma$; when $\delta \subseteq \sigma \subseteq \tau$, these homomorphisms satisfy $\eta^\tau_\delta = \eta^\sigma_\delta \circ \eta^\tau_\sigma$:  
 \beq  \label{Intertwining homoms diagram}  \xymatrixrowsep{2pc} \xymatrixcolsep{3pc}
\xymatrix{ G  \ar@<.5ex>[rr]^{\eta_\tau} \ar@<.5ex>[rrd]_{\eta_\sigma}
 & & G_\tau  \ar@{->}[d]^{\eta^\tau_\sigma } 
  \\
& & G_\sigma }
 \qquad \qquad \qquad  \xymatrixrowsep{2pc} \xymatrixcolsep{3pc}
\xymatrix{ G_\tau  \ar@<.5ex>[rr]^{\eta^\tau_\sigma} \ar@<.5ex>[rrd]_{\eta^\tau_\delta}
 & & G_\sigma \ar@{->}[d]^{\eta^\sigma_\delta } 
  \\
& & G_\delta }  \eeq 
When $\S$  is a  directed subset of $(\TG, \subseteq)$,  $\{ G_\tau; \ \eta^\tau_\sigma \text{ for } \sigma \subseteq \tau \text{ in } \S \}$ is thus a projective system of locally compact groups, and  we denote the associated   projective limit group by $G_\S$:   $$G_\S = \varprojlim_{\sigma \in {\cal S}} G_\sigma = \left\{ (s_\sigma)_{\sigma \in \S} \in \prod_{\sigma \in \S} G_\sigma:  \eta^\tau_\sigma(s_\tau) = s_\sigma \text{ for } \sigma \subseteq \tau \text{ in } \S \right\},$$  a closed topological subgroup of  $\prod_{\sigma \in \S} G_\sigma$.

\section{Preliminary results}

For $\tau$ in $\TG$,  $\eta_\tau: G \ra G_\tau$ is a continuous, dense-range homomorphism, so  $j_{\eta_\tau}: A(G_\tau) \ra B(G): v \mapsto v \circ \eta_\tau$ is a completely isometric isomorphism mapping $A(G_\tau)$ onto the Banach subalgebra $A_\tau(G) := A(G_\tau) \circ \eta_\tau$  of $B(G)$ \cite[Section 3]{Il-Spr1}. If $\lambda_\tau:= \lambda_{G_\tau} \circ \eta_\tau$, a continuous unitary representation of $G$ on $L^2(G_\tau)$, then $A_\tau(G) = A_{\lambda_{G_\tau}} \circ \eta_\tau$ is the Fourier space  $A_{\lambda_\tau}$ \cite[Propostion 2.10]{Ars}. Since the spectrum of $A(G_\tau)$ is $G_\tau$, we obtain the identification  \beq \label{Delta(Atau(G)Identitification} 
\Delta(A_\tau(G))= G_\tau \quad \text{via} \quad \l s, u\r_{G_\tau-A_\tau(G)} = \l s, \widehat{u} \r_{G_\tau-A(G_\tau)} = \widehat{u}(s),
\eeq
where $\widehat{\cdot}= j_{\eta_\tau}^{-1}: A_\tau(G) \ra A(G_\tau): u = v\circ \eta_\tau \mapsto \widehat{u} = v$.
If $\sigma, \tau $ belong to $\TG$, following \cite{Il-Spr1} we call $\sigma $ a \it quotient \rm of $\tau$ if $\sigma \subseteq \tau$ and $\eta^\sigma_\tau: G_\tau \ra G_\sigma$ is a proper map, meaning that pre-images of compact sets are compact; we will write $\sigma \preq \tau$ to indicate that $\sigma $ is a quotient of $\tau$. Obviously:

\blem \rm  The relation $\preq$ is a partial ordering of $\TG$. 
\elem

The   non-quotient topologies   in $\TG$, introduced in \cite{Il-Spr1} and denoted therein by $\TnqG$,  can thus be defined as the  maximal elements \rm in   $(\TG, \preq)$.  For $\tau \in \TG$, observe that $\tau$ is precompact --- meaning that $G_\tau$ is compact ---  if and only if   $(G_\tau, \eta_\tau)$  is a topological group compactification of $G$. If $\sigma, \tau \in \TG$ are precompact, notice that $\sigma \preq \tau$ exactly when $(G_\sigma, \eta_\sigma) \leq  (G_\tau, \eta_\tau)$ in the usual ordering of compactifications, e.g. see \cite{Ber-Jun-Mil}. Since the almost periodic compactification of $G$, $(G^{ap}, \eta_{ap})$, is the maximum  group compactification of $G$, $\tau_{ap} : = (\eta_{ap})^{-1} (\tau_{G_{ap}})$ is the maximum  precompact topology in $(\TG, \preq)$. Moreover, $\tau_{ap} \in \TnqG$: if $\tau_{ap} \preq \tau$, then $\tau_{ap} \subseteq \tau$ and $G_\tau = (\eta^\tau_{\tau_{ap}})^{-1}(G_{\tau_{ap}})$ is compact, so $\tau$ is precompact and therefore contained in $\tau_{ap}$.

   With the following proposition, we collect a number of facts about  $(\TG, \preq)$ and $\TnqG$, some new and some that are already found in \cite{Il-Spr1, Il-Spr2} as indicated in the proof.  Almost all of these statements are employed later in the paper, but some are included because we expect them to be  useful and of independent interest within the theory of spine algebras. 

\bp \label{Main preq proposition}  \rm Let $\tau_1, \tau_2, \sigma, \sigma', \tau, \tau', \delta $ be topologies in  $\TG$. The following statements hold:
\begin{enumerate}
\item $A_{\tau_1}(G) A_{\tau_2} (G) \subseteq  \overline{\text{span}(A_{\tau_1}(G) A_{\tau_2} (G))}^{\| \cdot \|_{B(G)}} = A_{\tau_1\vee \tau_2}(G)$.
\item  $\sigma \preq \tau$ if and only if $A_\sigma(G) \subseteq A_\tau(G)$; thus, $\sigma$ is uniquely determined by $A_\sigma(G)$. 
\item If $\sigma \preq \sigma'$ and $\tau \preq \tau'$, then $\sigma \vee \tau \preq \sigma'\vee \tau'$.

\item If $\delta \subseteq \sigma \subseteq \tau$, then $\delta \preq \tau$ if and only if $\delta \preq \sigma$ and $\sigma \preq \tau$.

\item If $\S$ is a directed subset of $(\TG, \subseteq)$, then   $\eta_\S: G \ra G_\S = \varprojlim_{\sigma \in \S} G_\sigma : s \mapsto (\eta_\sigma(s))_{\sigma \in \S}$ is a continuous dense-range homomorphism. 

\item If $\S$ is a directed subset of $(\TG, \preq)$, then   $G_\S = \varprojlim_{\sigma \in {\cal S}} G_\sigma$ is a locally compact group. Letting  $\tau_\S : = \eta_\S^{-1}(\tau_{G_\S})$, the weak topology induced by $\eta_\S$,  $\tau_\S$ belongs to  $\TG$ and  has representation $(G_\S, \eta_\S)$, $\tau_\S$ is the least upper bound of $\S$ in $(\TG, \subseteq)$, and   $\sigma \preq \tau_\S$ for each $\sigma$ in $\S$. 

\item  There is a unique topology $\taunq $ in $\TnqG$ such that $\tau \preq \taunq$. 
\item The set of topologies $Q_\tau := \{ \sigma \in \TG: \tau \preq \sigma\}$ is a $\vee$-subsemilattice of $\TG$ and a  directed subset of $(\TG, \preq)$. Moreover, $\taunq = \tau_{Q_\tau} =  \vee Q_\tau = \max Q_\tau$, where $\vee$ and $\max$ are taken in the lattice of  topologies on $G$ (with respect to containment).  
\item  If $\sigma\preq \tau$, then $\signq = \taunq$.
\item $(\sigma \vee \tau)_{\text{nq}} = (\signq \vee \taunq)_{\text{nq}}$; thus $(\TnqG, \widetilde{\vee})$, where $\sigma \widetilde{\vee} \tau:= (\sigma\vee \tau)_\nq$  is a join semilattice.
\item If $\sigma \subseteq \tau$, then $\signq \subseteq \taunq$. Hence, in $(\TnqG, \subseteq)$,  $\signq = \min \{ \tau\in \TnqG: \sigma \subseteq \tau\}$.

\end{enumerate}

\ep

\begin{proof} 1. This is \cite[Proposition 3.1]{Il-Spr1}, however it may not be true that $\lambda_{\tau_1} \otimes \lambda_{\tau_2}$ and $\lambda_{\tau_1 \vee \tau_2}$ are unitarily equivalent. (For example, if $G$ is a finite group and we take $\tau_1 = \tau_2 = \tau_G$, then $G_{\tau_1 \vee \tau_2} \cong G$, so   $\lambda_{\tau_1 \vee \tau_2}$ and $\lambda_{\tau_1} \otimes \lambda_{\tau_2} $ are, respectively, representations on  the non-isomorphic finite dimensional spaces $\ell^2(G)$ and $\ell^2( G \times G)$.)   We  therefore provide a  modified argument here, where we employ standard properties of Fourier algebras found, for example, in Sections 2.4 and 2.6 of \cite{Kan-Lau}:  

For $i=1,2$, let $u_i \in A_{\tau_i}(G)$ and take $u_i' \in A(G_{\tau_i})$ such that $u_i = j_{\eta_{\tau_i}}(u_i')$. Let $\lambda_i =  \lambda_{G_{\tau_i}}$ so $\lambda_{\tau_i} = \lambda_i \circ \eta_{\tau_i}$ and the outer tensor product representation $\{ \lambda_1 \times \lambda_2, L^2(G_{\tau_1}) \otimes L^2(G_{\tau_2})\}$ is unitarily equivalent to the left regular representation of $G_{\tau_1} \times G_{\tau_2}$. Taking $\xi_i, \ \eta_i \in L^2(G_{\tau_i})$ such that $u_i' = \phi^{\lambda_i}_{\xi_i, \eta_i}$, observe that $u_i = \phi^{\lambda_i}_{\xi_i, \eta_i} \circ \eta_{\tau_i} = \phi^{\lambda_{\tau_i}}_{\xi_i, \eta_i}$. Letting $\ov{\xi}_0 = \xi_1 \otimes \xi_2$, $\ov{\eta}_0 = \eta_1 \otimes \eta_2 \in L^2(G_{\tau_1}) \times L^2(G_{\tau_2})$, $$\phi^{\lambda_1 \times \lambda_2}_{\ov{\xi}_0, \ov{\eta}_0} \in A(G_{\tau_1} \times G_{\tau_2}), \quad \text{so} \quad  \phi^{\lambda_1 \times \lambda_2}_{\ov{\xi}_0, \ov{\eta}_0} \big{\vert}_{G_{\tau_1 \vee \tau_2}} \in A(G_{\tau_1 \vee \tau_2}).$$ But 
\beq \label{Product of Ataus Eqn} u_1u_2 =  j_{\eta_{\tau_1 \vee \tau_2}} \left(  \phi^{\lambda_1 \times \lambda_2}_{\ov{\xi}_0, \ov{\eta}_0} \Big{\vert}_{G_{\tau_1 \vee \tau_2}}  \right), \eeq
so $u_1 u_2 \in A_{\tau_1 \vee \tau_2}(G)$. Let $v \in A_{\tau_1 \vee \tau_2}(G)$, and take $w \in A(G_{\tau_1 \vee \tau_2})$ such that $v = j_{\eta_{\tau_1 \vee \tau_2}} w$, and $w' \in A(G_{\tau_1}  \times G_{\tau_2})$ such that $w'\big{\vert}_{G_{\tau_1 \vee \tau_2}} = w$. Take $\ov{\xi}, \ov{\eta}$ in $L^2(G_{\tau_1}) \otimes L^2(G_{\tau_2})$ such that $w' = \phi^{\lambda_1 \times \lambda_2}_{\ov{\xi}, \ov{\eta}}$, say $\ov{\xi} = \sum \ov{\xi}_n$, $\ov{\eta} = \sum \ov{\eta}_n$ where $\ov{\xi}_n = \xi^1_n \otimes \xi_n^2, \ \ov{\eta}_n = \eta_n^1 \otimes \eta_n^2$ are elementary tensors in $L^2(G_{\tau_1}) \otimes L^2(G_{\tau_2})$. 
Then $$\ds w' = \lim_N \sum_{n,m=1}^N \phi_{\ov{\xi}_n, \ov{\eta}_m}^{\lambda_1 \times \lambda_2} \quad \text{and therefore} \quad w = \lim_N \sum_{n,m=1}^N \left(\phi_{\ov{\xi}_n, \ov{\eta}_m}^{\lambda_1 \times \lambda_2} \Big{\vert}_{G_{\tau_1 \vee \tau_2}} \right).$$ Hence, $\ds v = j_{\eta_{\tau_1 \vee \tau_2}}w =  \lim_N \sum_{n,m=1}^N  j_{\eta_{\tau_1 \vee \tau_2}} \left(\phi_{\ov{\xi}_n, \ov{\eta}_m}^{\lambda_1 \times \lambda_2} \Big{\vert}_{G_{\tau_1 \vee \tau_2}} \right)$,  which 
 belongs to the closed linear span of $A_{\tau_1}(G) A_{\tau_2}(G)$  by (\ref{Product of Ataus Eqn}). 

\smallskip 

\noindent 2. This is an immediate consequence of Lemma 3.2 and Proposition 3.3 of  \cite{Il-Spr1}.

\smallskip 

\noindent 3.  This, which contains  \cite[Lemma 1.2(i)]{Il-Spr2}, follows quickly from parts 1 and 2:  We have $A_\sigma(G)A_\tau(G) \subseteq A_{\sigma'}(G)A_{\tau'}(G) \subseteq A_{\sigma'\vee\tau'}(G)$, so $A_{\sigma \vee \tau}(G) \subseteq A_{\sigma'\vee\tau'}(G)$; hence   $\sigma \vee \tau \preq \sigma' \vee \tau'$. 

\smallskip

\noindent 4. Suppose that $\delta \subseteq \sigma \subseteq \tau$ and $\delta \preq \tau$. That $\sigma \preq \tau$ is \cite[Lemma 1.2(ii)]{Il-Spr2}. Here we give a short proof that this is true, and also show that $\delta \preq \sigma$:  Let $C$ be a compact subset of $G_\sigma$. Since $C \subseteq (\eta_\delta^\sigma)^{-1} (\eta_\delta^\sigma(C))$, 
$$(\eta_\sigma^\tau)^{-1} (C)\subseteq (\eta_\sigma^\tau)^{-1} ( (\eta_\delta^\sigma)^{-1} (\eta_\delta^\sigma(C))) = (\eta_\delta^\sigma \circ \eta_\sigma^\tau)^{-1} (\eta_\delta^\sigma(C))
= (\eta_\delta^\tau)^{-1} (\eta_\delta^\sigma(C)).$$
Since $\eta_\delta^\sigma(C)$ is a compact subset of $G_\delta$ and $\eta^\tau_\delta$ is a proper map, we conclude that the closed set $(\eta_\sigma^\tau)^{-1} (C)$ is compact. Hence, $\sigma\preq\tau$.  To see that $\delta \preq \sigma$, take $K$ to be a compact subset of $G_\delta$. We just showed that $\eta^\tau_\sigma$ is a proper (therefore closed) dense-range map, hence a surjection, so 
$$(\eta^\sigma_\delta)^{-1}(K) = \eta^\tau_\sigma((\eta^\tau_\sigma)^{-1}((\eta^\sigma_\delta)^{-1}(K))) = \eta^\tau_\sigma((\eta^\sigma_\delta\circ \eta^\tau_\sigma)^{-1}(K)) = \eta^\tau_\sigma ((\eta^\tau_\delta)^{-1}(K)). $$
But $\eta^\tau_\delta$ is proper and $\eta^\tau_\sigma$ is continuous, so $(\eta^\sigma_\delta)^{-1}(K)$ is compact. Thus, $\delta \preq \sigma$.

\smallskip

\noindent 5. From (\ref{Intertwining homoms diagram}), $\eta_\S$ maps into $G_\S$ and it is clear that $\eta_\S$ is a continuous homomorphism.  Let $\F$ be a finite subset of $\S$, and consider the basic open subset $W = \bigcap_{\tau \in \F} \pi_\tau^{-1}(U_\tau)$ of $G_\S$,  where $U_\tau$ is open in $G_\tau$ and $\pi_\tau: G_\S \ra G_\tau$ is the $\tau^{\text{th}}$-projection map.  To show  that $\eta_\S$ has dense range in $G_\S$, we assume that $W$ is nonempty and show that $W \cap \eta_\S(G) \neq \emptyset$: Take $x = (x_\sigma)_{\sigma \in \S}$ in $W$ and take $\tau_0$ in $\S$ such that $\tau \subseteq \tau_0$ for each $\tau$  in $\F$.  For $\tau \in \F$,  $\eta^{\tau_0}_\tau(x_{\tau_0}) = x_\tau \in U_\tau$, so $x_{\tau_0} \in U_0:=\bigcap_{\tau \in \F} (\eta^{\tau_0}_\tau)^{-1}(U_\tau)$, an open subset of $G_{\tau_0}$. Since $\eta_{\tau_0}(G)$ is dense in $G_{\tau_0}$, we can choose $s \in G$ such that $\eta_{\tau_0}(s) \in U_0$. For $\tau \in \F$, $\pi_\tau(\eta_\S(s)) = \eta_\tau(s) = \eta^{\tau_0}_\tau(\eta_{\tau_0}(s)) \in G_\tau$, so $\eta_\S(s) \in W$, as needed. 

\smallskip 

\noindent 6. It follows from part  4 that $\S$ is a directed subset of $(\TG, \preq)$ if and only if it is a directed subset of $(\TG, \subseteq)$  with the property  that $\sigma \subseteq \tau$ in $\S$ implies that $\sigma\preq \tau$; so, $G_\S$ is a locally compact group  by \cite[Proposition 2.1]{Il-Spr1}, $\tau_\S \in \TG$,  and $(G_\S, \tau_\S)$ is a representation of $\tau_\S$ by statement 5.   Let $\sigma \in \S$, $V \in \sigma$, so $V = \eta_\sigma^{-1}(U)$ for some open subset $U$ of $G_\sigma$. Since $\eta_\sigma = \pi_\sigma \circ \eta_\S$, $V = \eta_\S^{-1}(\pi_\sigma^{-1}(U)) \in \tau_\S$. Hence, $\tau_\S$ is an upper bound of $\S$. Let $\tau$ be any upper bound of $\S$ in $\TG$. From (\ref{Intertwining homoms diagram}), $\eta^\tau_\S(x):= (\eta^\tau_\sigma(x))_{\sigma \in \S}$ maps $G_\tau$ into $G_\S$ and for $\sigma \in \S$, $\eta^\tau_\sigma = \pi_\sigma \circ \eta^\tau_\S$, so $\eta^\tau_\S$ is continuous, and $\eta^\tau_\S \circ \eta_\tau = \eta_\S$. So,  if $V \in \tau_\S$, say $V = \eta_\S^{-1}(U)$ where  $U$ is open in $G_\S$, then 
$V= \eta_\tau^{-1}((\eta^\tau_\S)^{-1}(U)) \in \tau$. Hence, $\tau_\S$ is the least upper bound of $\S$ in $(\TG, \subseteq)$.  Let $\sigma_0 \in \S$. To see that $\sigma_0 \preq \tau_\S$, first note that the projection homomorphism $\pi_{\sigma_0}: G_\S \ra G_{\sigma_0}$ satisfies $\pi_{\sigma_0} \circ \eta_\S = \eta_{\sigma_0}$, so $\eta^{\tau_\S}_{\sigma_0} = \pi_{\sigma_0}$. Let  $C$ to be a compact subset of $G_{\sigma_0}$. For each $\sigma \in \S$, choose $\tau_\sigma \in \S$ such that $\sigma, \sigma_0 \preq \tau_\sigma$, and observe that $C_\sigma := \eta^{\tau_\sigma}_\sigma ((\eta^{\tau_\sigma}_{\sigma_0})^{-1}(C))$ is a compact subset of $G_\sigma$. Since $G_\S$ is closed in $\prod_{\sigma \in \S} G_\sigma$, $D = (\prod_{\sigma \in \S} C_\sigma) \cap G_\S$ is a compact subset of $G_\S$.  We see that $\eta^{\tau_\S}_{\sigma_0}$ is proper by noting that the closed subset $(\eta^{\tau_\S}_{\sigma_0})^{-1}(C)$ of $G_\S$ is contained in $D$: Let $x = (x_\sigma)_{\sigma \in \S} \in (\eta^{\tau_\S}_{\sigma_0})^{-1}(C)$ and let $\sigma \in \S$. Then
$\eta^{\tau_\sigma}_{\sigma_0}(x_{\tau_\sigma}) = x_{\sigma_0} = \eta^{\tau_\S}_{\sigma_0}(x)  \in C$, so $x_{\tau_\sigma} \in (\eta^{\tau_\sigma}_{\sigma_0})^{-1}(C)$; hence $x_\sigma = \eta^{\tau_\sigma}_\sigma(x_{\tau_\sigma}) \in C_\sigma$. Thus, $x \in D$.  

\smallskip 

\noindent 7. This is \cite[Theorem 1.1(i)]{Il-Spr2}. However, the existence of $\taunq$ is also an easy consequence of  Zorn's lemma and part 6. Uniqueness follows from an  application of  property 3 (with $\sigma = \tau$), then 4.

\smallskip 

\noindent 8. The first sentence is a consequence of statement 3. By part 6,    $\tau_{Q_\tau} \in Q_\tau$ and $\sigma \preq \tau_{Q_\tau}$ for each $\sigma \in Q_\tau$. Hence, $\tau_{Q_\tau} = \max Q_\tau$. From statement 7, $\taunq \in Q_\tau$, so $\taunq = \tau_{Q_\tau}$.  Statement 8 is also implicit in  the proof  of \cite[Theorem 1.1(i)]{Il-Spr2}.

\smallskip 

\noindent 9. Suppose that $\sigma \preq \tau$. Then $Q_\tau \subseteq Q_\sigma$, so $\taunq = \max Q_\tau \subseteq \max Q_\sigma = \signq$.  Since  $\sigma \preq \signq$ and $\sigma \preq \tau$, $\sigma \preq \signq \vee \tau$ by part 3,  and $\sigma \subseteq \tau \subseteq \signq \vee \tau$, so $\tau \preq \signq \vee \tau$ by part 4.  Hence, $\signq \subseteq \signq \vee \tau \subseteq \max Q_\tau = \taunq$.

\smallskip 

\noindent 10.  Part 3 gives $\sigma \vee \tau \preq \signq \vee \taunq$, so this follows from part 9.

\smallskip 

\noindent 11.  If $\sigma \subseteq \tau$, then part 10 gives $\signq \subseteq (\signq \vee \taunq)_\nq = (\sigma \vee \tau)_\nq = \taunq$.  \end{proof} 


 Proposition \ref{Disjoint topologies Prop}, below, will be used in the proof of Theorem \ref{GSA MainThm2 Uniqueness of B}, but is also of independent interest. In the  lemma, $A(G:K)$ denotes the algebra of functions in $A(G)$ that are constant on the cosets of a compact normal  subgroup  $K$ of $G$.  Letting $q_K: G \ra G/K$ be the quotient homomorphism,  $j_{q_K}: A(G/K) \ra A(G:K)$ is a surjective isometric algebra isomorphism by \cite[Proposition 3.25]{Eym}. Moreover, $$P_K: A(G) \ra A(G:K) \text{ defined by } P_Ku(x) = \int_K u(xk) \, dm_K(k) \ (x \in G),$$ 
where $m_K$  is normalized Haar measure on $K$, is a projection of $A(G)$ onto $A(G:K)$ \cite[Corollary 3.4]{For}.

\blem  \label{Constant on cosets Lemma} \rm If $K$ and $L$ are compact normal subgroups of $G$, then $A(G:K) \cap A(G:L) \neq \{0\}$. 
\elem

\begin{proof} Take $u$ in $A(G)$ such that $u \equiv 1$ on the the compact subset $LK$ of $G$. Then $w:= P_L(P_K u))$ belongs to $A(G:K) \cap A(G:L)$ and $w(e_G) =1$. 
 \end{proof}

\bp \label{Disjoint topologies Prop} \rm Let $\sigma, \tau$ be topologies in $\TG$.  The following statements are equivalent:  
\bi 
\item[(i)] $A_\sigma(G) \cap A_\tau(G) \neq \{0\}$; 
\item[(ii)] there is a topology $\delta$  in $\TG$ such that $\delta \preq \sigma, \tau$; 
\item[(iii)] $\signq = \taunq$; 
\item[(iv)] there is a topology $\gamma$  in $\TG$ such that $\sigma, \tau \preq \gamma$.   
\ei 
\ep 

\begin{proof}  If (i) holds, then $A_\delta(G) = A_\sigma(G) \cap A_\tau(G)$ for some $\delta \preq \sigma, \tau$ by \cite[Lemma 3.2]{Il-Spr1}. By Proposition \ref{Main preq proposition}.10, statement (ii)  implies (iii) and (iv) follows trivially from (iii).  To see that (iv) implies (i), we begin by assuming that $\sigma \preq \gamma$ and show that $A_\sigma(G) = j_{\eta_\gamma}(A(G_\gamma : K_\sigma))$ where $K_\sigma = \ker \eta_\sigma^\gamma$. Indeed, $K_\sigma$ is a compact normal subgroup of $G_\gamma$ and $\tilde{ \eta}_\sigma^\gamma: G_\gamma/K_\sigma \ra G_\sigma$ is a topological group isomorphism \cite[p. 278]{Il-Spr1} so, letting $q_\sigma: G_\gamma \ra G_\gamma/K_\sigma$ be the quotient map, 
$$j_{ \eta_\sigma^\gamma} = j_{q_\sigma} \circ j_{\tilde{ \eta}_\sigma^\gamma}: A(G_\sigma) \ra A(G_\gamma/K_\sigma) \ra A(G_\gamma: K_\sigma)$$is an isometric algebra isomorphism of $A(G_\sigma)$ onto $A(G_\gamma: K_\sigma)$. Since $\eta_\sigma =  \eta_\sigma^\gamma \circ \eta_\gamma$, 
$$A_\sigma (G) = j_{\eta_\sigma}(A(G_\sigma)) = j_{\eta_\gamma}(j_{ \eta_\sigma^\gamma}(A(G_\sigma))) = j_{\eta_\gamma}(A(G_\gamma:K_\sigma)),$$ as claimed.  If we also assume that $\tau \preq \gamma$, then $A_\tau(G) = j_{\eta_\gamma}(A(G_\gamma: K_\tau))$ where $K_\tau = \ker  \eta_\tau^\gamma$, so statement (i) follows from Lemma \ref{Constant on cosets Lemma} because $j_{\eta_\gamma}: A(G_\gamma) \ra A_\gamma(G)$ is injective.   
   \end{proof} 

\br  \rm The above proof contains a refinement of Proposition \ref{Main preq proposition}.2 that seems worth recording: If $\sigma \preq \gamma$, then $A_\sigma(G) = j_{\eta_\gamma}(A(G_\gamma : K_\sigma)) \subseteq A_\gamma(G)$ where $K_\sigma = \ker \eta_\sigma^\gamma$.

\er

\section{Generalized spine algebras and their spectra}

For a set $\D$ of pre-locally compact topologies on $G$, i.e., for $\D \subseteq \TG$, let    $$\ADSG = \overline{\sum_{\tau \in \D} A_\tau(G)},$$ 
the norm-closed linear span of $\{ A_\tau(G): \tau \in \D\}$ in $B(G)$.  We will say that a subset $\B$ of $\TG$ is  \it disjoint \rm if the associated representations $\{\lambda_\tau: \tau \in \B\}$ are disjoint, equivalently if $A_\sigma (G) \cap A_\tau(G)=\{0\}$ for $\sigma \neq \tau$ in $\B$. In this case, 
$$A_\B^*(G) = \ell^1-\bigoplus_{\tau \in \B} \,  A_\tau(G) = A_{\lambda_\B}$$ 
where $\{ \lambda_\B, \H_\B\} = \{ \bigoplus_{\tau \in \B} \, \lambda_\tau, \bigoplus_{\tau \in \B} \,  L^2(G_\tau)\}$ and each $u\in \ABSG$ has a unique representation \beq \label{Unique rep as sum over disjoint tops Eqn} u =  \sum_{\tau \in \B} u_\tau \ \  \text{ where}  \ u_\tau \in A_\tau(G) \text{ and, moreover, } \  \|u\| = \sum_{\tau \in \B} \| u_\tau\| \eeq   by \cite[3.12 and 3.13]{Ars}.  For example, $\TnqG$ is disjoint by \cite[Lemma 3.4]{Il-Spr1} and the full spine subalgebra of $B(G)$ 
is $$A^*(G) = A^*_{\TG}(G) = A^*_{\TnqG}(G).$$

\bd \rm A subspace  $\A$ of $B(G)$ is a \it generalized spine algebra (GSA) \rm if  $\A$ is an algebra and  $\A = \ADSG$ for some subset $\D$ of $\TG$.    \ed

\subsection{Direct sum decomposition of a GSA} 

Let $\A$ be a GSA. Let $$\D_\A := \{ \sigma \in \TG: A_\sigma(G) \subseteq \A\} \quad \text{and} \quad \TnqA := \{ \tau \in \TnqG: A_\tau(G) \cap \A \neq \{0\}\}.$$ For $\tau$ in $\TnqA$, let $\tau_\A$ in $\TG$ be the unique locally precompact topology on $G$ such that 
$$A_{\tau_\A}(G) = A_\tau(G) \cap \A$$
\cite[Lemma 3.2]{Il-Spr1} and let 
$$\B_\A := \{ \tau_\A: \tau \in \TnqA\}.$$
For $\sigma \in \D_\A$, $\signq \in \TnqA$ and we let 
$$\sigma^\sharp = (\signq)_\A, \text{ an element of }  \B_\A.$$ 

For a GSA $\A$, $\D_\A$ and $\B_\A$ take the roles played by  $\TG$ and $\TnqG$ for the full spine algebra $A^*(G)$.  

\bp   \label{SigmaSharpProp}  \rm  Let $\A$ be a GSA. 

\begin{enumerate}  \item For $\sigma \in \D_\A$,  $\sigma \preq \sigma^\sharp \preq \signq$  and 
$\signq = (\sigma^\sharp)_\nq$.

\item For $\sigma, \tau \in \D_\A$,  $\sigma^\sharp = \tau^\sharp$ if and only if $\signq = \taunq$. 
\smallskip 

\item  $\B_\A =\{ \tau_\A : \tau \in \TnqA\}$ is a disjoint subset of $\D_\A$. 
\item If $\sigma \in \D_\A$, then $\sigma^\sharp$ is the unique topology in $\B_\A$ satisfying $\sigma \preq \sigma^\sharp$ and $\tau^\sharp = \tau$ for $\tau \in \B_\A$.

\item $\B_\A$ is the set of maximal elements in the partially ordered set $(\D_\A, \preq)$. 

\end{enumerate}  
\ep

\begin{proof}   For $\sigma \in \D_\A$,  $A_\sigma(G) \subseteq A_\signq(G) \cap \A = A_{\sigma^\sharp}(G) \subseteq A_{\signq}(G)$, so $\sigma \preq \sigma^\sharp \preq \signq$ and  $\signq = (\sigma^\sharp)_\nq$ by parts 2 and 9 of  Proposition \ref{Main preq proposition}.   If $\sigma, \tau \in \D_\A$ and $\signq =\taunq$, then $A_{\sigma^\sharp}(G) = A_{\tau^\sharp}(G)$, so $\sigma^\sharp = \tau^\sharp$, again by Proposition \ref{Main preq proposition}.2. 
For $\tau$ in $\TnqA$, $\tau_\A \preq \tau$ so $(\tau_\A)_\nq = \tau$, and we conclude that $\B_\A$ is disjoint because $\TnqA$ is so. We have established statements 1--3 and statement 4 is a consequence of parts 1, 3 and Proposition  \ref{Main preq proposition}.2. Finally, take $\sigma \in \B_\A$, $\tau \in \D_\A$ and suppose that $\sigma \preq \tau$. Then $\signq = \taunq$ by Proposition \ref{Main preq proposition}.9, so parts 1, 2 and 4 yield $\tau^\sharp = \sigma^\sharp = \sigma \preq \tau \preq \tau^\sharp$; hence $\sigma = \tau$. 
\end{proof}

Observe that if $\D$ is a join subsemilattice of $(\TG, \vee)$ --- we write $\D \leq (\TG, \vee)$ when this is the case ---  then $\ADSG$ is a GSA because $A_\sigma(G) A_\tau(G) \subseteq A_{\sigma \vee \tau}(G)$.

\blem \label{A = AD for semilattice D Lemma}  \rm Let $\A$ be a GSA. Then $\D_\A \leq (\TG, \vee)$, $\A = A_{\D_\A}^*(G)$ and if $\D$ is any join subsemilattice of $(\TG, \vee)$ such that $\A = A_\D^*(G)$, then $\D \leq \D_\A$. 
\elem

\begin{proof}  If ${\cal E}$ is a subset of $\TG$ such that $\A = A_{\cal E}^*(G)$, then ${\cal E} \subseteq \D_\A$, so $\A = A_{\D_\A}^*(G)$; by Proposition \ref{Main preq proposition}.1,  $\D_\A \leq (\TG, \vee)$. \end{proof}


\bp \label{GSAProp1}  \rm 1. If $\A$ is a GSA and $\sigma, \tau$ in $\D_\A$, let $$\sigma \vee_\A \tau: = (\sigma \vee \tau)^\sharp.$$ 
 Then  $\sigma \vee_\A\tau = (\sigma \vee \tau)^\sharp =  (\sigma^\sharp \vee \tau^\sharp)^\sharp = \sigma^\sharp \vee_\A \tau^\sharp$.

\medskip 

\noindent 2.  For a subset $\A$ of $B(G)$, the following statements are equivalent: 
\bi  \item[(i)] $\A$ is a GSA; 
\item[(ii)] there is a join subsemilattice $\D$ of $(\TG, \vee)$ such that $\A = \ADSG$; 
\item[(iii)] there is a  disjoint set of  locally precompact topologies $\B$ such that $\A = \ABSG$ and for every $\sigma, \tau $ in $\B$ there is a (necessarily unique) topology $\gamma_{\sigma, \tau}$ in $\B$ such that $\sigma \vee \tau \preq \gamma_{\sigma, \tau}$.  
\ei  
\noindent 3. Suppose that $\A= A_\D^*(G)$ is a GSA, where $\D$ is a join subsemilattice of  $(\TG, \vee)$. Then:
\bi \item[(a)]  $\D$ and $\B_\D := \{ \sigma^\sharp: \sigma \in \D\}$ are contained in $\D_\A$, and $\B_\D$ is a disjoint set of locally precompact topologies on $G$ such that statement 2 (iii) is satisfied with $\gamma_{\sigma, \tau} = \sigma \vee_\A \tau$.
\item[(b)] $(\B_\D, \vee_\A)$ is a join semilattice and  $$(\D, \vee) \ra (\B_\D, \vee_\A) : \sigma \mapsto \sigma^\sharp$$ is a join semilattice epimorphism. 
\ei 

\ep

\begin{proof}  We will employ various parts of Proposition \ref{SigmaSharpProp}.  

\smallskip  

\noindent 1. Let $\sigma, \tau \in \D_\A$. By Proposition \ref{Main preq proposition}.10, 
$$ (\sigma \vee \tau)_\nq = (\signq \vee \taunq)_\nq = ((\sigma^\sharp)_\nq \vee (\tau^\sharp)_\nq)_\nq = (\sigma^\sharp \vee \tau^\sharp)_\nq,$$ 
so $\sigma \vee_\A\tau = (\sigma \vee \tau)^\sharp = ((\sigma \vee \tau)_\nq)_\A = ((\sigma^\sharp \vee \tau^\sharp)_\nq)_\A  = (\sigma^\sharp \vee \tau^\sharp)^\sharp=  \sigma^\sharp \vee_\A \tau^\sharp.$  

\smallskip 

\noindent 2.  By Lemma  \ref{A = AD for semilattice D Lemma} and the comment preceding it, statements (i) and (ii) are equivalent,  and it is obvious that (iii) implies (i). Suppose that $\A = \ADSG$, where $\D \leq (\TG, \vee)$.  Note that  $\B_\D$ is disjoint since it  is contained in $\B_\A$ and, since $A_\sigma(G) \subseteq A_{\sigma^\sharp}(G) \subseteq \A = \ADSG$ ($\sigma \in \D$),  $\A = A_{\B_\D}^*(G)$. Let $\sigma, \tau \in \B_\D$, say $\sigma = \sigma_0^\sharp$ and $\tau = \tau_0^\sharp$ where $\sigma_0, \tau_0 \in \D$. From part 1,   $\sigma \vee_\A \tau  =  (\sigma_0 \vee \tau_0)^\sharp  \in \B_\D$ and $\sigma \vee \tau \preq \sigma \vee_\A \tau$.  If $\B$ is now taken to be any set of topologies satisfying  (iii), and $\sigma, \tau, \gamma_1, \gamma_2$ in $\B$ satisfy $\sigma \vee \tau \preq \gamma_1, \gamma_2$, then $\{0\} \neq A_\sigma(G) A_\tau(G) \subseteq A_{\sigma \vee \tau}(G) \subseteq A_{\gamma_1}(G) \cap A_{\gamma_2}(G)$;  so, $\gamma_1 = \gamma_2 $ because $\B$ is disjoint. Thus, $\gamma_{\sigma, \tau}$ is unique and $\sigma \vee_\A \tau = \gamma_{\sigma, \tau}$ when $\B = \B_\D$.

\smallskip 

\noindent 3. Statement 3(a) has been justified and 3(b) follows from part 1. For example, as part of showing that $(\B_\D, \vee_\A)$ is a join semilattice, take $\gamma, \sigma, \tau \in \D$ and observe that 
$$ \gamma^\sharp \vee_\A (\sigma^\sharp \vee_\A\tau^\sharp )=  \gamma^\sharp \vee_\A (\sigma \vee\tau)^\sharp = (\gamma  \vee (\sigma  \vee \tau))^\sharp = ((\gamma  \vee \sigma)  \vee \tau)^\sharp   = (\gamma^\sharp \vee_\A \sigma^\sharp) \vee_\A\tau^\sharp. \qedhere$$ 
\end{proof} 


\bt  \label{GSA MainThm2 Uniqueness of B}   \rm Let $\A$ be a GSA. 
Then  $\B = \B_\A$ is the  unique    disjoint set of locally precompact topologies on $G$ such that 
$$\A = A_{\B}^*(G) = \ell^1-\bigoplus_{\tau  \in \B} A_\tau(G).$$  Moreover, $(\B_\A, \vee_\A)$ is a join semilattice such that for each $\sigma, \tau \in \B_\A$,  $\sigma \vee_\A \tau$ is the unique topology in $\B_\A$ satisfying $$A_\sigma(G) A_\tau (G) \subseteq A_{\sigma \vee_\A \tau}(G)  \quad  \text{(equivalently, } \sigma \vee \tau \preq \sigma \vee_\A \tau).$$  The map $(\D_\A, \vee)\ra (\B_\A,\vee_\A): \sigma \mapsto \sigma^\sharp$ is a join semilattice epimorphism with kernel $$\{(\sigma, \tau) \in \D_\A \times \D_\A: A_\sigma(G) \cap A_\tau(G) \neq \{0\}\},$$ and  $$\kappa: (\B_\A, \vee_\A) \mapsto (\TnqA, \widetilde{\vee} ): \sigma \mapsto \signq$$ is a join semilattice isomorphism with inverse map $\tau \mapsto \tau_\A$. 
\et  

\begin{proof}  Take $\D\leq (\TG, \vee)$ such that $\A = A_\D^*(G)$. Then $\B_\D \subseteq \B_\A$ and $\A = A_{\B_\D}^*(G)$ by Proposition \ref{GSAProp1}.3, so $\A = A_{\B_\A}^*(G)$.  Let $\B$ be any disjoint subset of $\TG$ such that $\A = A_\B^*(G)$.   Then 
\beq  \ell^1-\bigoplus_{\sigma \in \B} A_\sigma(G) = A_\B^*(G) & = & \A \label{GSA MainThm2 Eq1} \\
& = & A_{\B_\A}^*(G) =  \ell^1-\bigoplus_{\tau \in \B_\A} A_\tau(G) \label{GSA MainThm2 Eq2} \eeq
and for each $\sigma \in \B$, $\sigma^\sharp \in \B_\A$ and $A_\sigma(G) \subseteq A_{\sigma^\sharp}(G)$. Hence, the unique representation of any $u$ in $\A$ with respect to $\B$ via (\ref{GSA MainThm2 Eq1}), $\sum_{\tau \in \B} u_\tau$,  is also the unique representation of $u$ with respect to $\B_\A$ via (\ref{GSA MainThm2 Eq2}), so  $A_\sigma(G) = A_{\sigma^\sharp}(G)$, equivalently $\sigma = \sigma^\sharp$, for each $\sigma$ in $\B$. (Indeed, if there were some $u$ in $A_{\sigma^\sharp}(G) \bs A_\sigma(G)$, $u$ would have two distinct representations with respect to $\B_\A$.) Hence, $\B \subseteq \B_\A$.  If there were some $\tau \in \B_\A \bs  \B$, any $u \in A_\tau(G)$ would have two representations $\sum u_\sigma$ taken with respect to $\B_\A$, so $\B = \B_\A$.  It now follows from Proposition \ref{GSAProp1} that $\B_\A = \B_\D$, $(\B_\A, \vee_\A)$ is a semilattice and $\sigma\mapsto \sigma^\sharp$ is a semilattice epimorphism of $\D_\A$ onto $\B_\A$. By Proposition \ref{SigmaSharpProp}.2,  $\sigma^\sharp = \tau^\sharp$ if and only if $\signq = \taunq$, so the description of the kernel of this map is a consequence of Proposition  \ref{Disjoint topologies Prop}.   We only need to check that $\kappa$ is  a join isomorphism. Certainly, $\kappa$ maps $\B_\A$ into $\TnqA$ and it  is one-to-one  by Proposition \ref{SigmaSharpProp}. For $\sigma, \tau$ in $\B_\A$, 
$$\kappa(\sigma \vee_\A \tau) = ((\sigma \vee \tau)^\sharp)_\nq = (\sigma \vee \tau)_\nq = (\signq \vee \taunq)_\nq = \kappa(\sigma) \widetilde{\vee} \kappa(\tau),$$
where we have used Propositions \ref{SigmaSharpProp}.1 and \ref{Main preq proposition}.10.   For $\tau \in \TnqA$, $\tau_\A \in \B_\A$ and $\tau_\A \preq \tau$, so $(\tau_\A)_\nq = \tau$; hence,   $\kappa$ is also surjective and $\kappa^{-1}(\tau) = \tau_\A$. 
\end{proof}

\bc \label{GSA Main Thm Cor1}  \rm  Suppose that $(\B, \vee)$ is a disjoint subsemilattice of $\TG$. Then $\A = A_\B^*(G) = \ell^1-\bigoplus_{\tau \in \B} A_\tau(G) $ is a GSA, $\B_\A  =\B$, and $\sigma \vee \tau = \sigma \vee_\A \tau$ for each $\sigma, \tau \in \B$.  \ec 

\bc \rm  \label{GSA Main Thm Cor} Let $\A$ be a GSA,  $\sigma$ and $\tau$ in $\B_\A$. Then $\sigma \vee_\A \tau = \tau$ if and only if $\sigma \subseteq \tau$; i.e., $\subseteq$ is the partial ordering on $\B_\A$ induced by $\vee_\A$. So, with respect to $\subseteq$,  $\sigma \vee_\A \tau$ is the least upper bound in $\B_\A$ of $\sigma$ and $\tau$. 

 \ec 
 
 \begin{proof}   
 Clearly, if $\sigma \vee_\A\tau = \tau$, then $\sigma \subseteq \tau$ and if $\sigma \subseteq \tau$, then $\sigma \vee_\A \tau = (\sigma \vee \tau)^\sharp = \tau^\sharp = \tau$ by Proposition \ref{SigmaSharpProp}.4. 
 \end{proof} 
 
 \br  \rm 1.  In Section 4, we will show how to construct many new examples of GSAs. 
The spine and reduced spine algebras, $A^*(G)$ and  $A^*_0(G) = A^*_{{\cal T}_0(G)}(G)$ where a locally precompact topology $\tau$ belongs to $ {\cal T}_0(G) $ if it is Hausdorff,  are GSAs that were studied in \cite{Il-Spr1}.  A particularly accessible GSA considered in \cite{Sto21}
is $A_{\mathbb{F}}(G) = A(G) \oplus_1 A(G^{ap}) \circ \eta_{ap} =  A_\B^*(G)$, where $\B= \{ \tau_{ap}, \tau_G\}$. In \cite{Il-Spr1} it is shown that $A^*(G) =  A_{\mathbb{F}}(G)$ for some groups, including $\R$, $\mathbb{Z}$, the $p$-adics $\mathbb{Q}_p$, and any minimally weakly almost periodic group such as $SL_2(\R)$ or a Euclidean motion group. 

\smallskip 

\noindent 2. GSAs, in the case of abelian groups, are examples of certain balanced $L$-subalgebras $\mathfrak{N}$ of commutative convolution measure algebras $\mathfrak{M}$ introduced by Joseph Taylor. These subalgebras are characterized by the condition that the spine of ${\mathfrak{N}}$, as defined in \cite{Tay},  coincides with the intersection  of $\mathfrak{N}$ with the spine of ${\mathfrak{M}}$. Taylor reduced the invertibility problem in an abelian measure algebra to the same problem in a  much less complicated balanced $L$-subalgebra \cite[Chapter 9]{Tay}. 
\er

\subsection{Gelfand spectrum of a GSA}

Letting $\A$ be a GSA, we will identify its spectrum $\Delta(\A)$. Our approach, though slightly different from that found of \cite{Il-Spr1},  is based on ideas from Section 4.1 therein, which in turn were influenced by ideas in  \cite{Inoue}.   Going forward, we will usually write $\B$ in place of $\B_\A$, so  
$\A = A_\B^*(G) = A_{\lambda_\B}$ where $\lambda_\B = \oplus_{\tau \in \B} \lambda_\tau$.   
Writing $VN_\tau = VN_{\lambda_\tau}$, 
$\A^* = VN_{\lambda_\B} = \ell^\infty-\bigoplus_{\tau \in \B} VN_\tau$ through the pairing
 \beq \label{GSA Dual Pairing Eqn} \l x, u \r = \sum_{\tau \in \B} \l x_\tau, u_\tau\r \quad \text{for } x = (x_\tau)_{\tau \in \B} \in VN_{\lambda_\B}, \ u = \sum_{\tau \in \B} u_\tau \in A_\B^*(G) = \ell^1-\bigoplus_{\tau \in \B} A_\tau(G) \eeq
  \cite[(3.13) Cor.]{Ars}. If $t =(t_\tau)_{\tau \in \B} \in \Delta(\A) \subseteq VN_{\lambda_\B}$, then for each $\tau$ in $\B$, $t_\tau = t\large{\vert}_{A_\tau(G) } \in \Delta(A_\tau(G)) \cup \{ 0 \} = G_\tau \cup \{ 0 \};$ hence $$\Delta(A_\B^*(G)) \subseteq \Pi_{\tau \in \B} (G_\tau \cup \{0\}) \subseteq VN_{\lambda_\B}.$$  For $\S \subseteq \B$, we may write 
  \beq \label{Spectrum Notation Eqn} t = (t_\tau)_{\tau \in \S} \in \Pi_{\tau \in \B} (G_\tau \cup \{0\}) \text{ when } t_\tau = 0 \text{ for } \tau \in \B \bs \S.\eeq
  
 A nonempty subset $\S$ of $\B$ is \it hereditary and directed \rm if $\S$ is closed downward in $\B$ with resect to containment  and is a join subsemilattice of $(\B, \vee_\A)$ (equivalently, a directed subset of $(\B, \subseteq)$) cf. \cite{Il-Spr1}. Let $\text{HD}(\B)$ be the collection of all hereditary and directed subsets of $\B$. For $\sigma \in \B$, $\S_\sigma:= \{\tau \in \B: \tau \subseteq \sigma\} $  is a \it principal \rm element of $\text{HD}(\B)$ and we let $\PHD(\B) =\{ \S_\sigma: \sigma \in \B\}$. For $\S$ in $\text{HD}(\B)$, $G_\S = \varprojlim_{\sigma \in \S} G_\sigma$ is a topological group contained in $\Pi_{\tau \in \B} (G_\tau\cup\{0\}) \subseteq VN_{\lambda_\B}$, using our convention (\ref{Spectrum Notation Eqn}). Thus, for $x \in G_\S$, $x \in VN_{\lambda_\B} = (A_\B^*(G))^*$ through the pairing (\ref{GSA Dual Pairing Eqn}). Let $$G_\B^* = \bigcupdot_{\S \in \text{HD}(\B)} G_\S.$$
As noted in Section 2, we are using $\widehat{\cdot}$ to denote the inverse of the isometric algebra isomorphism $j_{\eta_\tau}: A(G_\tau) \mapsto A_\tau(G)$ and the  identification $\Delta(A_\tau(G)) = G_\tau$ is described by (\ref{Delta(Atau(G)Identitification}).

 \bp  \label{Gelfand Spectrum Main Prop}  \rm  Let $\A = A_{\B}^*(G)$ be a GSA, where $\B = \B_\A$. Then $\ds \Delta(\A) = G_\B^*$.
 \ep

 \begin{proof}   Let $\sigma, \tau \in \B$, $u \in A_\sigma(G)$, $v \in A_\tau(G)$. Then $uv \in A_{\sigma \vee_\A \tau}(G)$ and the argument used to establish (4.3) of \cite{Il-Spr1} shows that for any $t \in G_{\sigma \vee_\A \tau}$, 
 \beq \label{Ilie-Spronk Eqn}  \widehat{uv}(t) = \widehat{u}(\eta^{\sigma\vee_\A\tau}_\sigma(t))  \widehat{v} (\eta^{\sigma\vee_\A\tau}_\tau(t)).
 \eeq

 Suppose that $s = (s_\sigma)_{\sigma \in \S} \in G_\S$. As noted above, $s \in VN_{\lambda_\B} = \A^*$, so it suffices to show that $s$ is multiplicative. Take $u, v$ as above. If $\sigma, \tau \in \S$, then $\sigma\vee_\A\tau \in \S$, so  from (\ref{Delta(Atau(G)Identitification}), (\ref{GSA Dual Pairing Eqn}), (\ref{Ilie-Spronk Eqn})  and the definition of $G_\S$, we obtain 
 \beqs  \l s, uv\r &  = & \l s_{\sigma \vee_\A \tau}, uv \r = \widehat{uv}(s_{\sigma\vee_\A\tau}) = \widehat{u}(\eta^{\sigma\vee_\A\tau}_\sigma(s_{\sigma\vee_\A\tau}))  \widehat{v} (\eta^{\sigma\vee_\A\tau}_\tau(s_{\sigma\vee_\A\tau}))\\
 & = & \widehat{u}(s_\sigma )  \widehat{v} (s_\tau)  =   \l s_\sigma, u \r \l s_\tau, v\r =  \l s, u \r \l s, v\r.
 \eeqs 
 If $\sigma \notin \S$ or $\tau \notin \S$, then $\sigma \vee_\A \tau \notin \S$ and we again obtain $\l s, uv \r  = \l s, u \r \l s, v\r (=0)$. 
 Since $  \A = \ell^1-\bigoplus_{\tau \in \B} A_\tau (G) = \ov{\sum_{\tau \in \B} A_\tau(G)} $, we conclude that $s \in \Delta (\A)$.

 Conversely, suppose that $s  \in \Delta(\A)$. As noted above, $s= (s_\tau)_{\tau \in \B} \in \Pi_{\tau \in \B} (G_\tau \cup \{0\})$ and, letting $\S = \{ \tau \in \B: s_\tau \neq 0\}$, we can write $s = (s_\tau)_{\tau \in \S}$. To see that $\S \in \text{HD}(\B)$ and $s \in G_\S$, first take $\sigma, \tau \in \S$ and $u \in A_\sigma(G)$, $v \in A_\tau(G)$ such that $\l s_\sigma, u \r = \l s_\tau, v \r  = 1$. Then $uv \in A_{\sigma \vee_\A \tau}(G)$, so $$\l s_{\sigma \vee_\A\tau}, uv \r = \l s , uv \r = \l s, u \r \l s, v \r = \l s_\sigma, u \r \l s_\tau, v \r = 1;$$ hence $\sigma \vee_\A \tau \in \S$. Suppose now that $\sigma \in \S$ and $\gamma \in \B$ with $\gamma \subseteq \sigma$. Again, take $u \in A_\sigma(G)$ such that $\l s, u \r = \l s_\sigma, u \r =1$, and note that $s_\gamma \in G_\gamma \cup \{ 0\} \subseteq VN_\gamma$, $\eta^\sigma_\gamma(s_\sigma)  \in  G_\gamma$. For any $v \in A_\gamma(G)$, $uv \in A_\sigma(G)$ and (\ref{Ilie-Spronk Eqn}) yields, 
 \beqs  \l s_\gamma, v \r &= & \l s, u \r \l s, v \r = \l s, uv \r = \l s_\sigma, uv \r = \widehat{uv} (s_\sigma) \\ 
 & = & \widehat{u}(\eta^\sigma_\sigma(s_\sigma)) \widehat{v} ( \eta^\sigma_\gamma(s_\sigma)) = \l s_\sigma, u \r \l \eta^\sigma_\gamma(s_\sigma), v\r = \l \eta^\sigma_\gamma(s_\sigma), v\r.
 \eeqs
 Hence, $s_\gamma = \eta^\sigma_\gamma(s_\sigma) \in G_\gamma$ (and $\gamma \in \S$), as needed. 
 \end{proof} 
 
For us, a semigroup  $$S = \bigcupdot_{z \in Z} G_z$$ is a \it semilattice of disjoint semigroups \rm if  
\bi \item $(Z, \cdot, \leq)$ is a (lower) semilattice, i.e., a commutative semigroup of idempotents with ordering $w\leq z$ if $wz = w$; and 
\item the disjoint semigroups $G_z$ are graded over $Z$, in the sense that $G_z G_w \subseteq G_{zw}$. 
\ei 
We will say that $S$ has \it compatible central identities \rm if the semigroups $G_z$ have  identities $1_z$ that are  central in $S$ and satisfy $1_z1_w = 1_{zw}$ for $z, w$  in $Z$. 

\br \rm When $S =  \bigcupdot_{z \in Z} G_z$ is a semilattice of groups with compatible central identities, then the semilattice $Z$ and the semilattice of idempotents of $S$,  $E(S)$, are isomorphic via  $z\mapsto 1_z$; in this case, $E(S)$ can thus be used as an index set in place of $Z$.   \er

Let $\A = A_{\B}^*(G)$ be a GSA, where $\B = \B_\A$. By \cite[Corollary 1.6]{Sto21} and its proof,  $\Delta(\A) \cup \{0\} = G_\B^* \cup \{0 \}$ is a multiplicative $*$-subsemigroup of $\A^* = VN_{\lambda_{\B}}$ and when $\A$ contains the identity function $1_G$, $\Delta(\A) = G_\B^*$ is itself a $*$-semigroup. Using our convention (\ref{Spectrum Notation Eqn}), $G_\S = \{0\}$ when $\S = \emptyset$, so $$\Delta(\A) \cup \{0\} =   \bigcupdot_{\S \in \HD(\B) \cup \{\emptyset\}}  G_\S.$$  
The product and involution on $\Delta(\A) \cup \{0\} $, inherited from $VN_{\lambda_{\B}} = \ell^\infty-\bigoplus_{\tau \in \B} VN_\tau$, is given coordinatewise:  for $s= (s_\sigma)_{\sigma \in \S}, \ t = (t_\tau)_{\tau \in \cal{S}'} \in \Delta(\A) \cup \{0\}$, 
\beq   \label{Product in Delta(A) Eqn} st = (s_\sigma t_\sigma)_{\sigma \in \S \cap \cal{S}'} \quad \text{and} \quad  s^* = (s_\sigma^*)_{\sigma \in \S}  = (s_\sigma^{-1})_{\sigma \in \S}.
  \eeq     
For $\sigma \in \TG$, let $e_\sigma = e_{G_\sigma}$, the identity of $G_\sigma$, so $e_\S = (e_\sigma)_{\sigma \in \S}$ is the identity of the group  $G_\S$.   Observe that $\HD(\B) \cup \{ \emptyset\}$ is closed under arbitrary intersections.  

\bp \label{Delta(A)SemigpProp}  \rm Let $\A = A_{\B}^*(G)$ be a GSA, where $\B= \B_\A$.  \bi \item[1.]
The spectrum of $\A$, $\Delta(\A) = G_\B^*=  \bigcupdot_{\S \in \HD(\B)}  G_\S$, is a $*$-semigroup under the operations (\ref{Product in Delta(A) Eqn}) if and only if $$\S \cap \S' \neq \emptyset \quad \text{ for each }  \S, \ \S' \in \HD(\B).$$ When this condition is satisfied, 
 $G_\B^*= \bigcupdot_{\S \in \HD(\B)}  G_\S$ is a semilattice of groups, graded over the lower semilattice $(\HD(\B), \cap, \subseteq)$, with compatible central identities $e_\S = (e_\sigma)_{\sigma \in \S}$. Also,   $(\HD(\B), \subseteq)$  has maximum element $\B$ and is a complete upper semilattice with $\vee {\mathfrak Q}$ given by $$\vee {\mathfrak Q} = \bigcap\{{\cal T} \in \HD(\B): \S \subseteq {\cal T} \text{ for each } \S \in {\mathfrak Q}\}$$ for any nonempty subset ${\mathfrak Q}$ of   $\HD(\B)$.   The map $\S \mapsto e_\S$ is a semilattice isomorphism of $(\HD(\B), \cap)$ onto $E(\Delta(\A))$, the semilattice of idempotents (= projections) in $\Delta(\A)$.
  
 \item[2.] If $(\B, \subseteq)$ has a minimum element or if $1_G \in \A$, then $\Delta(\A) = G_\B^*$ is a semilattice of groups, graded over the (complete when $\min \B$ exists) lower semilattice $(\HD(\B), \cap, \subseteq)$. 
 
 \item[3.] It is always true that $\Delta(\A) \cup \{0\} = G_\B^*\cup\{0\} =   \bigcupdot_{\S \in \HD(\B) \cup \{\emptyset\}}  G_\S$ is a semilattice of groups, graded over the complete lower semilattice $(\HD(\B)\cup\{\emptyset\}, \cap, \subseteq)$, with compatible central identities $e_\S = (e_\sigma)_{\sigma \in \S}$. 
\ei 
\ep 

\begin{proof} This is not difficult to check. For example, observe that if the intersection of any two sets in $\HD(\B)$ is nonempty, then $\HD(\B)$ is closed under finite intersections, so $(\HD(\B), \cap, \subseteq)$ is a semilattice. In this case, if $s,t \in \Delta(\A)$, say  $s= (s_\sigma)_{\sigma \in \S} \in G_\S$ and  $t = (t_\tau)_{\tau \in \S'} \in G_{\S'}$, then 
$st =    (s_\sigma t_\sigma)_{\sigma \in \S \cap \S'} \neq 0$ and if $\sigma \subseteq \tau$ in $\S \cap \S'$, then 
$\eta^\tau_\sigma(s_\tau t_\tau) =  \eta^\tau_\sigma(s_\tau) \eta^\tau_\sigma( t_\tau) = s_\sigma t_\sigma,$ so $st \in G_{\S \cap \S'}$. Thus $G_\S G_{\S'} \subseteq G_{\S \cap \S'} $ and $e_\S e_{\S'}= e_{\S \cap \S'}$. 
\end{proof} 

Let $\A$ be the full spine algebra over $G$, $A^*(G)$. Notice  that  $1_G \in A_{\tau_{ap}}(G) \subseteq  \A$. Also,  $\tau_{ap}$ is the minimum element of $\TnqG = \B_\A$: if $\tau \in \TnqG $ and $v \in A_\tau(G)$, then $v = 1_G \, v \in A_{\tau_{ap} \vee_\A \tau} (G) \cap A_\tau(G)$ and $\TnqG$ is disjoint, so $\tau_{ap} \vee_\A  \tau = \tau$; hence $\tau_{ap} \subseteq \tau$.  A special case of  Propositions \ref{Gelfand Spectrum Main Prop} and \ref{Delta(A)SemigpProp}  is, thus,  the following corollary, which is contained in  \cite[Section 4.1]{Il-Spr1}. 

\bc  \rm  The spectrum, $\Delta(A^*(G))$,  of the full spine algebra $A^*(G)$  is $$ G^*= \bigcupdot_{\S \in \HD(\TnqG)} G_\S,$$ which is a semilattice of groups, graded over the complete lower semilattice $(\HD(\TnqG), \cap, \subseteq)$.  \ec  

                 
\section{Omnibus example} 

For any upper subsemilattice $\cal{E}$  of $({\cal P}(\N), \subseteq, \cup)$, we will show how to build a  GSA, $\A = A_\B^*(G)$, for which the disjoint upper semilattice of  locally precompact topologies $(\B=\B_\A, \vee_\A)$ is isomorphic to  $(\cal{E}, \cup)$ and $\vee_\A = \vee$, the ordinary join in $\TG$. Given a GSA, $\A$, the interesting semigroup structure of $\Delta(\A) = \bigcupdot_{\S \in \HD(\B)}  G_\S$ is, as it pertains to the goals of this paper, primarily determined by the lower semilattice $( \HD(\B), \cap)$.  As we will see, our procedure  allows us to build a class of subalgebras of Fourier--Stieltjes algebras with  a wide variety---one for each subsemilattice $\cal E$ of ${\cal P}(\N)$---of complicated, yet often rather tractable, spectra. For instance, with Example  \ref{P(N) Example} we will show how to build $\A$ so that  $( \HD(\B), \cap)$ is a lower-semilattice isomorphic with $({\cal P}(\N), \cap)$.  

 If we begin with an infinite abelian group $H$, with the following omnibus construction we will produce a disjoint semilattice  $\{ \tau_A: A \subseteq \N\}$  of Hausdorff  locally precompact topologies on $G = H^\N$; see Theorem \ref{MainOmnibusThm} and Remark \ref{NonHausdorffOmnibus Remark}.1. A simpler method for producing a semilattice  $\{ \sigma_A: A \subseteq \N\}$ of  locally precompact topologies on $G = H^\N$  with similar properties, but such that the topologies $\sigma_A$ are never Hausdorff, is described in Remark \ref{NonHausdorffOmnibus Remark}.2.  

    For clarity, in this section we use  $(G^{ap}, \eta^G_{ap})$ for the almost periodic compactification of $G$ and   $\tau_{ap}^G$ is the corresponding topology in $\TG$.

\blem \rm  \label{Product of AP Cpctns Lemma} Let $(G, \tau_G)$ and $(H, \tau_H)$ be locally compact groups. 
\bi  \item[1.] Let $\phi:(G, \tau_G) \ra (H, \tau_H)$ be a topological group isomorphism.  Then there is a topological group isomorphism $\widetilde{\phi}: G^{ap} \ra H^{ap}$ such that $\widetilde{\phi} \circ \eta^G_{ap} = \eta^H_{ap} \circ \phi$. Moreover,  $\phi: (G, \tau_{ap}^G) \ra (H, \tau_{ap}^H)$ is also a topological group isomorphism. 

\item[2.] The topological group compactifications  $((G\times H)^{ap}, \eta^{G \times H}_{ap})$ and $(G^{ap} \times H^{ap}, \eta^{G}_{ap} \times \eta^H_{ap})$ of $G \times H$ are equivalent. Thus,  $(G^{ap} \times H^{ap}, \eta^{G}_{ap} \times \eta^H_{ap})$  is a representation of the topology  $\tau^{G \times H}_{ap}$ in ${\cal T}( G \times H)$ and  $(g_\alpha, h_\alpha) \ra (g, h)$ in $\tau^{G \times H}_{ap}$ if and only if $\eta^G_{ap}(g_\alpha) \ra  \eta^G_{ap}(g)$ in $G^{ap}$ and $\eta^H_{ap}(h_\alpha) \ra  \eta^H_{ap}(h)$ in $H^{ap}$. 
\ei   
\elem  

\begin{proof} 1. Since $(H^{ap}, \eta^H_{ap} \circ \phi)$ and $(G^{ap}, \eta^G_{ap} \circ \phi^{-1})$ are topological group compactifications of $G$ and $H$ respectively, the universal property of an almost periodic compactification yields continuous surjective homomorphisms  $\widetilde{\phi}: G^{ap} \ra H^{ap}$ and $\widetilde{\psi}: H^{ap} \ra G^{ap}$ satisfying $\widetilde{\phi} \circ \eta^G_{ap} = \eta^H_{ap} \circ \phi$ and $\widetilde{\psi} \circ \eta^H_{ap} = \eta^G_{ap} \circ \phi^{-1}$. For $g \in G$, $\widetilde{\psi} \circ \widetilde{\phi}( \eta^G_{ap}(g)) = \widetilde{\psi}(\eta^H_{ap}(\phi(g))) = \eta^G_{ap}(\phi^{-1}(\phi(g))) = \eta^G_{ap}(g).$  Similarly, $\widetilde{\phi} \circ \widetilde{\psi}$ is the identity map on the dense subset $\eta^H_{ap}(H)$ of $H^{ap}$, so $\widetilde{\phi}$ is a topological group isomorphism (with inverse $\widetilde{\psi}$).  From (\ref{NetCvgceTauEqn}), the second statement in part 1 follows because  $g_\alpha \ra g$ in $\tau^G_{ap}$ if and only if $\eta^G_{ap}(g_\alpha) \ra \eta^G_{ap}(g)$ in $G^{ap}$, which holds exactly when $\eta^H_{ap}(\phi(g_\alpha)) = \widetilde{\phi}(\eta^G_{ap}(g_\alpha)) \ra \widetilde{\phi}(\eta^G_{ap}(g))  = \eta^H_{ap}(\phi(g))$; equivalently, $\phi(g_\alpha) \ra \phi(g)$ in $\tau^H_{ap}$. 

\smallskip 

\noindent 2.  We observe that the topological group compactification $(G^{ap} \times H^{ap}, \eta^G_{ap} \times \eta^H_{ap})$ of $G\times H$ is  universal among all topological group compactifications of $G \times H$. (This must be  known but we have no reference.) To this end, let $((G\times H)^\alpha, \alpha)$ be a group compactification of $G \times H$. Then $(G^{\alpha_G}, \alpha_G)$ and $(H^{\alpha_H}, \alpha_H)$ are group compactifications of $G$ and $H$ respectively, where $\alpha_G, \alpha_H: G, H \ra (G \times H)^\alpha$ are defined by $\alpha_G(g) = \alpha(g, e_H)$, $\alpha_H(h) = \alpha(e_G, h)$ and $G^{\alpha_G}$, $H^{\alpha_H}$ are, respectively, the closures in $(G \times H)^\alpha$ of $\alpha_G(G)$ and $\alpha_H(H)$. The universal properties of $G^{ap}$ and $H^{ap}$  yield continuous (group homomorphisms) $\phi_G: G^{ap} \ra G^{\alpha_G}$ and $\phi_H: H^{ap} \ra H^{\alpha_H}$ satisfying $\phi_G \circ \eta^G_{ap} = \alpha_G$ and $\phi_H \circ \eta^H_{ap} = \alpha_H$. The map $\phi: G^{ap} \times H^{ap} \ra (G \times H)^\alpha: (s,t) \mapsto \phi_G(s) \phi_H(t)$ is  continuous and satisfies $\phi \circ (\eta^G_{ap} \times \eta^H_{ap}) = \alpha$, as needed; $\phi$ is automatically a surjective homomorphism.   Since $(G^{ap} \times H^{ap}, \eta^{G}_{ap} \times \eta^H_{ap})$ is equivalent to the almost periodic compactification of $G \times H$, it is  a representation of $\tau^{G \times H}_{ap}$.  The final statement follows from (\ref{NetCvgceTauEqn}). 
 \end{proof} 

 Throughout the remainder of this section,  $H$ is a fixed infinite discrete  group  and  $G= H^\N$ with the discrete topology.  Let $A, B$ be subsets of $\N$. If $A \subseteq B$, let 
 $$\pi_A : H^B \ra H^A : (x_i)_{i \in B} \mapsto  (x_i)_{i \in A}$$  and put $$\phi_A: G \ra H^A \times H^{\N\bs A}: x \mapsto (\pi_A(x), \pi_{\N \bs A}(x)).$$
When $A$ and $B$ are disjoint, let 
 $$\lambda_{A, B} : H^A \times H^B \ra H^{A \cup B}: ((x_i)_{i \in A}, (x_i)_{i \in B}) \mapsto (x_i)_{i \in A \cup B},$$
  $$\phi_{A, B} :  H^{A \cup B} \ra  H^A \times H^B: x \mapsto (\pi_A(x), \pi_B(x)).$$
The group $H^A$ is given the discrete topology and we will write $\eta^A_{ap}$ in place of the AP-compactification homomorphism $\eta^{H^A}_{ap}: H^A \ra (H^A)^{ap}$, etc. Let 
$$G_A = H^A \times (H^{\N \bs A})^{ap}$$
(with $G_\N = H^\N = G$, $G_\emptyset = G^{ap}$) and put $$\eta_A = (\text{id}_A \times \eta_{ap}^{\N \bs A}) \circ \phi_A: G \ra H^A \times H^{\N \bs A} \ra G_A : x \mapsto (\pi_A(x) , \eta^{\N \bs A}_{ap}(\pi_{\N \bs A}(x))). $$ 
Then $\eta_A$ is a continuous dense-range homomorphism, so 
$$\tau_A:= \eta_A^{-1} (\tau_{G_A}) \in \TG; \quad \tau_\N = \tau_G, \quad  \tau_\emptyset = \tau^G_{ap}.$$

 \blem \rm \label{Cvge in tauA Lemma} 1.  For $A \subseteq \N$, $x_i \ra x$ in $(G, \tau_A)$ if and only if $\pi_A(x_i) \ra \pi_A(x)$ in $H^A$  and $\eta^{\N \bs A}_{ap}(\pi_{\N\bs A}(x_i)) \ra \eta^{\N \bs A}_{ap}(\pi_{\N \bs A}(x))$. 
 
 \smallskip  
 
 \noindent 2.  Let $E, F$ be disjoint subsets of $\N$, $(x_i)$, $(y_i)$, $(z_i)$ and $x, \ y, z$ nets and elements in $H^E$, $H^F$ and $H^{E\cup F}$, respectively.  Then: \bi\item[(a)] $\lim_i \eta^{E}_{ap}(x_i) = \eta^{E}_{ap}(x)$ in $(H^E)^{ap}$ and $\lim_i \eta^{F}_{ap}(y_i) = \eta^{E}_{ap}(y)$ in $(H^F)^{ap}$ if and only if $\lim_i\eta^{E \cup F}_{ap}(\lambda_{E, F}(x_i, y_i)) = \eta^{E \cup F}_{ap}(\lambda_{E,F}(x,y)$  in $(H^{E\cup F})^{ap}$. 
\item[(b)]  $\lim_i\eta^{E \cup F}_{ap}(z_i) = \eta^{E \cup F}_{ap}(z)$  in $(H^{E\cup F})^{ap}$ if and only if   $\lim_i\eta^{E}_{ap}(\pi_E(z_i)) = \eta^{E}_{ap}(\pi_E(z)) $ in $(H^{E})^{ap}$ and  $\lim_i\eta^{F}_{ap}(\pi_F(z_i)) = \eta^{F}_{ap}(\pi_F(z)) $ in $(H^{F})^{ap}$.
\ei 
 \elem
 
 \begin{proof}    Statement 1 is clear and statement 2(a) follows from Lemma \ref{Product of AP Cpctns Lemma} because $H^{A \cup B}= H^A \times H^B$. Statement 2(b) is equivalent to statement 2(a). 
 \end{proof} 
 
 Let $A$ and $B$ be subsets of $\N$  with $A \subseteq B$.  Let $\kappa_{A, B}$ be the topological group isomorphism defined through the commuting diagram 
 \beqs   \xymatrixrowsep{2pc} \xymatrixcolsep{2pc}
\xymatrix{ (H^{B\bs A})^{ap} \times (H^{\N \bs B})^{ap} \ar@<.5ex>[d]_{ \psi} \ar@<.5ex>[rr]^{ \quad \kappa_{A,B}}  & & (H^{\N \bs A})^{ap}
  \\
(H^{B\bs A} \times H^{\N \bs B})^{ap} \ar@<.5ex>[rru]_{ \ \ \ \widetilde{\lambda}_{B \bs A, \N \bs B}}   }   \eeqs
 where  $\psi$ is the compactification isomorphism (see Lemma \ref{Product of AP Cpctns Lemma}.2) and 
$\widetilde{\lambda}_{B \bs A, \N \bs B}$ is the topological group isomorphism  induced by the map $\lambda_{B\bs A, \N \bs B}$ via Lemma \ref{Product of AP Cpctns Lemma}.1. Observe that 
\beq  \label{kappaABEqn}
\kappa_{A,B} \circ  (\eta_{ap}^{B\bs A} \times \eta_{ap}^{\N \bs B})  = \widetilde{\lambda}_{B \bs A, \N \bs B} \circ \eta_{ap}^{H^{B\bs A} \times H^{\N \bs B}} = \eta_{ap}^{\N \bs A} \circ \lambda_{B \bs A, \N \bs B}. 
\eeq

 \bt  \label{MainOmnibusThm}  \rm  1. The map $${\cal P}(\N)  \ra \TG: A \mapsto \tau_A$$ 
 is injective and $\tau_{A\cup B} = \tau_A \vee \tau_B$, (i.e., $A \mapsto \tau_A$ is an upper semilattice monomorphism). 
 
 \smallskip 
 
 \no 2. For $A\subseteq B$, $\tau_A \subseteq \tau_B$ and  the canonical intertwining map $\eta^
  B_A := \eta^{\tau_B}_{\tau_A} $ factors as 
  \beqs   \xymatrixrowsep{2pc} \xymatrixcolsep{2pc}
\xymatrix{ G_B= H^B \times (H^{\N \bs B})^{ap} \ar@<.5ex>[d]^\cong_{\phi_{A,B\bs A} \times \text{id}} \ar@<.5ex>[rr]^{\eta_A^B}  & & G_A= H^A \times (H^{\N \bs A})^{ap}
  \\
H^A \times H^{B\bs A} \times (H^{\N \bs B})^{ap} \ar@<.5ex>[rr]^{\text{id}_A \times\eta_{ap}^{B\bs A} \times \text{id}}  & & H^A \times (H^{B\bs A})^{ap} \times (H^{\N \bs B})^{ap}     \ar@<.5ex>[u]^\cong_{\text{id}_A \times \kappa_{A,B}} }   \eeqs
 
 \smallskip 
 
 \no 3.  The set $\B = \{ \tau_A: A \subseteq \N \} $ is a disjoint subsemilattice of $(\TG, 
 \vee)$. 
 \et 
 
 \begin{proof}  1. Suppose that $A \subsetneq B$ and suppose that $x_i \ra x$ in $\tau_B$.  By Lemma \ref{Cvge in tauA Lemma}.1, $\pi_A(x_i) \ra \pi_A(x)$ in $H^A$ (with the discrete topology), $\pi_{B \bs A} (x_i) \ra \pi_{B\bs A}(x)$ in $H^{B \bs A}$ --- so, $\eta^{B \bs A}_{ap}(\pi_{B\bs A}(x_i)) \ra \eta^{B \bs A}_{ap}(\pi_{B\bs A}(x))$ --- and $\eta^{\N \bs B}_{ap}(\pi_{\N\bs B}(x_i)) \ra \eta^{\N \bs B}_{ap}(\pi_{\N \bs B}(x))$. Since $(\N \bs B) \cup (B\bs A) = \N \bs A$,  $\eta^{\N \bs A}_{ap}(\pi_{\N\bs A}(x_i)) \ra \eta^{\N \bs A}_{ap}(\pi_{\N \bs A}(x))$ by   Lemma \ref{Cvge in tauA Lemma}.2(b). We conclude that $x_i \ra x$ in $\tau_A$ and, therefore, $\tau_A \subseteq \tau_B$.   To see that $\tau_A \subsetneq \tau_B$,  we take $m \in B\bs A$ and observe that because $H$ is Hausdorff and noncompact, we can choose $h$ and a net $(h_i)$ in $H$ such that $\lim \eta_{ap}^H(h_i) = \eta_{ap}^H(h)$, but $\lim h_i \neq h$ in $H$. Letting  $\iota_m: H \ra G=H^\N$ be defined by $\iota_m(k)_n = k$ if $n=m$, $\iota_m(k)_n = e_H$ otherwise, one sees from Lemma \ref{Cvge in tauA Lemma} that $\lim \iota_m(h_i)  = \iota_m(h)$ in $\tau_A$  but $\lim \iota_m(h_i)  \neq \iota_m(h)$ in $\tau_B$. Thus, $\tau_A \subsetneq \tau_B$, as needed.

 Supposing now that $A$ and $B$ are any two subsets of $\N$, we know that $\tau_A, \tau_B \subseteq \tau_{A\cup B}$. Let $\sigma$ be a topology on $G$ that contains $\tau_A$ and $\tau_B$ and suppose that $x_i \ra x$ in $\sigma$. Then $x_i \ra x$ in $\tau_A$ and $\tau_B$, so part 1 of  Lemma \ref{Cvge in tauA Lemma}  clearly  yields $\pi_{A\cup B}(x_i) \ra \pi_{A \cup B}(x)$ in $H^{A\cup B}$; using part 1, then part 2(b), of Lemma \ref{Cvge in tauA Lemma},  we see that $\eta^{\N \bs (A\cup B)}_{ap}(\pi_{\N\bs (A\cup B)}(x_i)) \ra \eta^{\N \bs (A\cup B)}_{ap}(\pi_{\N \bs (A\cup B)}(x))$ because $\N \bs (A \cup B) \subseteq \N \bs A$. Hence,  $x_i \ra x$ in $\tau_{A\cup B}$ and therefore  $\tau_{A \cup B} \subseteq \sigma$. We conclude that $\tau_A \vee \tau_B = \tau_{A \cup B}$.  
 
 Finally, for part 1, suppose that $A \neq B$, say (without loss of generality)  $A \subsetneq A \cup B$. Then $\tau_A \subsetneq \tau_{A \cup B} = \tau_A \vee \tau_B$, so $\tau_A \neq \tau_B$.

 \smallskip 
 
 \noindent 2. Suppose that $A \subseteq B$, so that $\tau_A \subseteq \tau_B$. To see that $\eta^{\tau_B}_{\tau_A} = \eta^
  B_A$, as described in statement 2,  we will show that $\eta^B_A\circ \eta_B = \eta_A$.  For  $x$ in $G = H^\N$, 
\beqs  \eta^B_A \circ \eta_B(x) & = & (\text{id}_A \times  \kappa_{A,B}) \circ ( \text{id}_A  \times \eta^{B \bs A}_{ap} \times \text{id}) \circ (\phi_{A, B \bs A} \times \text{id})  (\pi_B(x), \eta_{ap}^{\N \bs B} (\pi_{\N \bs B}(x))) \\& = & (\text{id}_A \times  \kappa_{A,B})(\pi_A(x), (\eta^{B \bs A}_{ap} (\pi_{B \bs A}(x)), \eta^{\N \bs B}_{ap}(\pi_{\N \bs B}(x)))\\
& = & (\pi_A(x), \kappa_{A,B} \circ(\eta_{ap}^{B \bs A} \times \eta_{ap}^{\N \bs B}) (\pi_{B\bs A}(x), \pi_{\N \bs B}(x)) \\
& = & (\pi_A(x), \eta_{ap}^{\N \bs A} \circ \lambda_{B \bs A, \N \bs B}(\pi_{B\bs A}(x), \pi_{\N \bs B}(x)) \\
& = & (\pi_A(x), \eta_{ap}^{\N \bs A} (\pi_{\N \bs A}(x)))\\
& = & \eta_A(x),  
\eeqs   
where we have used equation (\ref{kappaABEqn}). 

\smallskip  

\noindent 3. Suppose that $A \subsetneq B$, so $\tau_A \subsetneq \tau_B$. To see that $\tau_A$ is not a quotient of $\tau_B$, i.e., to see  that $\eta^B_A$ is not a proper map, consider the compact subset $K = \{e_{H^A}\} \times (H^{\N \bs A})^{ap}$ of $G_A = H^A \times (H^{\N \bs A})^{ap}$.  Take $m \in B \bs A$, let $H_m = \{ h = (h_k)_{k \in B} : h_k = e_H \text{ for } k \neq m\}$ and put $L = H_m \times (H^{\N \bs B})^{ap}$. As $L$ is a closed, non-compact subset of $(\eta^B_A)^{-1}(K)$, we conclude that $\tau_A$ is not a quotient of $\tau_B$. 

  Finally, let $A,B$ be any two sets and suppose that $A_{\tau_A}(G) \cap A_{\tau_B}(G) \neq \{0\}$. By Proposition \ref{Disjoint topologies Prop}, there is topology $\gamma$ in $\TG$ such that $\tau_A, \tau_B \preq \gamma$. As $\tau_A, \tau_B \subseteq \tau_A \vee \tau_B \subseteq \gamma$, $\tau_A, \tau_B \preq \tau_A \vee \tau_B$ by Proposition \ref{Main preq proposition}.4. But $A, B \subseteq A\cup B$ and $\tau_A \vee \tau_B = \tau_{A \cup B}$, so we know from the previous case that $A = A \cup B = B$.  
 \end{proof}

\br \label{NonHausdorffOmnibus Remark}   \rm 1. In the above construction, if $H$ is taken to be an infinite \it abelian \rm discrete group, then $\eta_A: G \ra G_A$ is injective, so $\tau_A$ is Hausdorff for each subset $A$ of $\N$.

\smallskip 

\noindent 2. As before, let $H$ be an infinite discrete group, $G = H^\N$ with the discrete topology. With fewer notational hurdles, we can produce an upper semilattice of \it non-Hausdorff \rm topologies $\{ \sigma_A: A \subseteq \N\}$ in $\TG$ such that $A \mapsto \sigma_A : ({ \cal P}(\N), \cup, \subseteq) \ra (\TG, \vee, \subseteq)$ is a semilattice monomorphism as follows: 

For $A\subseteq \N$, give $H^A$ the discrete topology and let $\sigma_A$ be the weak topology  on $G$  induced by the projection homomorphism $\eta_A: G= H^\N \ra H^A$.  Then a nonempty set $U$ belongs to $\sigma_A$ if and only if $U = V \times H^{\N \bs A}$ for some nonempty subset $V$ of $H^A$ and $\lim x_\alpha = x$ in $\sigma_A$ exactly when  $\eta_A(x_\alpha) = \eta_A(x)$, eventually. From this,  one readily verifies that $\sigma_A \subsetneq \sigma_B$ for $A \subsetneq B$, $\sigma_{A \cup B} = \sigma_A \vee \sigma_B$, and $A \mapsto \sigma_A$ is one-to-one. For $A \subsetneq B$, $\eta^{\sigma_B}_{\sigma_A} = \sigma^B_A:H^B \ra H^A$ is the projection homomorphism, which has noncompact kernel and is therefore not a proper map; so, $\sigma_A$ is not a quotient of $\sigma_B$. The argument found in the last paragraph of the proof of Theorem \ref{MainOmnibusThm} now shows that $\{\sigma_A: A \subseteq \N\}$ is disjoint. Thus, going forward  one can use $\{ \sigma_A: A \subseteq \N\}$ as an alternative choice to $\{ \tau_A: A \subseteq \N\}$.  

\smallskip 

\noindent 3. For a subset $A$ of $\N$,  $A_{\sigma_A}(G) \cong A(H^A)$ and  $A_{\tau_A}(G) \cong A(H^A \times (H^{\N \bs A})^{ap} ) \cong A(H^A) \widehat{\otimes} A(H^{\N \bs A})^{ap} )$ through completely isometric algebra isomorphisms, where $\widehat{\otimes}$ denotes the operator space projective tensor product \cite{Eff-Rua}. 
\er
Let $\E$ be an upper subsemilattice of $({\cal P}(\N), \cup, \subseteq)$. By Theorem \ref{MainOmnibusThm}, $\B_\E= \{ \tau_A: A \in \E\}$ is a disjoint subsemilattice of $(\TG, \vee, \subseteq)$ that is join-isomorphic with $\E$ and, by Corollary \ref{GSA Main Thm Cor1}, $\A=A_{\B_\E}^*(G)$ is a GSA with $(\B_\A, \vee_\A)= (\B_\E, \vee)$. To simplify notation, we will write $A_\E^*(G)$ in place of  $A_{\B_\E}^*(G)$. Let $$\HD(\E) = \{ \S : \S \text{ is a subsemilattice of } (\E, \cup) \text{ and closed downwards in } (\E, \subseteq)\}.$$  Assume that $\S \cap \S'$ is nonempty for $\S, \S' \in \HD(\E)$, such as when $\E$ has a minimum element.  (If this is not the case, replace $\E$ with  $\E\cup\{\emptyset\}$.)  Observe that the lower semilattices $(\HD(\E), \cap)$ and $(\HD(\B_\E), \cap)$ are then isomorphic via $\S \mapsto \{ \tau_A: A \in \S\}$. 

\bc   \rm  \label{Omnibus Example Cor}  Let $\E$ be an upper subsemilattice of $({\cal P}(\N), \cup, \subseteq)$ such that $\S \cap \S'$ is nonempty for $\S, \S' \in \HD(\E)$. Then 
\bi \item[1.] $\ds A_\E^*(G) = \ell^1-\bigoplus_{A \in \E} A_{\tau_A}(G)$ is a GSA such that $A_{\tau_A}(G) A_{\tau_B}(G) \subseteq A_{\tau_{A\cup B}}(G)$;
\item[2.] $\Delta(A_\E^*(G)) = \bigcupdot_{ \S \in \HD(\E)} G_\S$, where 
$$G_\S = \varprojlim_{A \in \S} G_A = \{ (s_A)_{A \in \S}: \eta^B_A(s_B) = s_A \text{ for } A\subseteq B \text { in } \S\},$$
is a semilattice of disjoint groups with compatible central identities, graded over the lower semilattice $(\HD(\E), \cap)$; $(\HD(\E), \subseteq)$  is a complete upper semilattice with maximum element $\E$.
\ei 
\ec

By strategically choosing $\E$, one can often perfectly describe the lower semilattice $(\HD(\E), \cap)$, which, as noted above, from our perspective determines most of the interesting semigroup structure of $\Delta(A_\E^*(G)) = \bigcupdot_{\S \in \HD(\E)} G_\S$.  We begin with two examples where $(\HD(\E), \cap) \cong ({\cal P}([n]), \cap)$ and  $(\HD(\E), \cap) \cong ({\cal P}(\N), \cap)$.

 \bex \rm \label{Finite Boolean Alg Example} Let $n \in \N$, $[n]= \{1,2, \ldots, n\}$, and let $\E_n= \{ A \in {\cal P}(\N): A \subseteq [n]\}$, so $$(\E_n, \cup ) \cong ({\cal P}([n]), \cup).$$ In this case, if   $\S \in \HD(\E_n)$, then $A = \bigcup \S \in \S$, so $\S = \S_A = \{B : B \subseteq A\}$. Thus,  $\HD(\E_n) = \PHD(\E_n)$  and $\S_A \cap \S_B = \S_{A\cap B}$, so
$(\HD(\E_n), \cap) \cong ({\cal P}([n]), \cap)$. Hence,  $\Delta(A^*_{\E_n}(G)) = \bigcupdot_{A \subseteq [n]} G_{\S_A},$ where $G_{\S_A} G_{\S_B} \subseteq G_{\S_{A\cap B}}$ via $(s_C)_{C \subseteq A} (t_C)_{C \subseteq B} = (s_Ct_C)_{C \subseteq A\cap B}$. Since $A = \max \S_A$, 
$$G_A \cong G_{\S_A}  = \varprojlim_{B \subseteq A} G_B \quad \text{via} \quad s \mapsto (\eta^A_B(s))_{B \subseteq A},$$
$$\Delta(A^*_{\E_n}(G)) = \bigcupdot_{A \subseteq [n]} G_A \quad \text{via} \quad \l s, \sum_{B\subseteq [n]} u_B\r = \sum_{B \subseteq A} \l \eta^A_B(s), u_B\r \text{ for } s \in G_A,$$ and  $\Delta(A^*_{\E_n}(G)) $ is a  semilattice of disjoint groups with compatible central identities, graded over $({\cal P}([n]), \cap)$, with product 
\beq \label{Omnibus product in spectrum of G_A G_B Eqn}  G_A G_B \subseteq G_{A \cap B} \quad \text{given by} \quad s t= \eta^A_{A \cap B}(s) \eta^B_{A \cap B}(t)  \text{ for } s \in G_A, \ t \in G_B. \eeq
This example is illustrated in  Figure 1.  
\begin{figure} 
\begin{multicols}{3}

  \centering  \begin{tikzpicture}
  
    \matrix (A) [matrix of nodes, row sep=1.2cm]
    { 
        & $G_{\{1,2\}}$ \\  
	$G_{\{1\}}$ &  & $G_{\{2\}}$\\
        & $G_\emptyset = G^{ap}$ \\
    };
    \draw (A-1-2)--(A-2-1);
    \draw (A-1-2)--(A-2-3);
    \draw (A-2-1)--(A-3-2);
    \draw (A-2-3)--(A-3-2);
\end{tikzpicture} 

\columnbreak 

$$G_AG_B \subseteq G_{A\cap B}$$

\columnbreak 
\centering
\begin{tikzpicture}
    \matrix (A) [matrix of nodes, row sep=1.2cm]
    {
        $G_{\{1,2\}}$ & $ G_{\{1,3\}}$ & $G_{\{2,3\}}$ \\
        $G_{\{1\}}$ & $G_{\{2 \}}$ & $G_{\{3\}}$ \\
        & $G_\emptyset=G^{ap}$ \\
    };
    \path (A-1-1)--(A-1-2) node[above=1.2cm] (link) {$G_{\{1,2,3\}}$};
    
    \foreach \i in {1,...,3}
    \draw (link.south) -- (A-1-\i.north);
    \foreach \i/\j in {1/2, 3/2, 2/1, 1/1, 3/3, 2/3}
    \draw (A-1-\i.south)--(A-2-\j.north);
    \foreach \i/\j in {1/2, 2/2, 3/2}
    \draw (A-2-\i.south)--(A-3-\j.north);
\end{tikzpicture}
\end{multicols}
\caption{$\ds \Delta(A^*_{\E_n}(G)) = \bigcupdot_{A \subseteq [n]} G_A$ for $n = 2, 3$}\label{Figure1} 
\end{figure} 
 \eex
 
 \br \rm When $\E$ is a subsemilattice of $({\cal P}([n]), \cup, \subseteq)$, like in Example \ref{Finite Boolean Alg Example}, it should be clear that one can replace $G=H^\N$ with $G= H^n$.  
 \er

 \bex  \label{P(N) Example}  \rm  Let $\F= \{ A \subseteq \N: A \text{ is finite}\}$, an upper subsemilattice of $({\cal P}(\N), \cup, \subseteq)$. Observe that 
 $$\HD(\F) = \{ \S_B: B \subseteq \N\} \ \ \text{where} \ \ \S_B = \{ A \in \F:  A \subseteq B\}  \text{ for } B \subseteq \N.$$
 Indeed, if $\S \in \HD(\F)$, $\S = \S_B$ where $B = \bigcup \S$: taking $A$ in $\S_B$, so that $A$ is a \it finite \rm subset of $B = \bigcup \S$, for each $a \in A$ there is a set $F_a \in \S$ such that $a \in F_a$; hence, $A \subseteq \bigcup_{a \in A} F_a \in \S$ and therefore $A \in \S$ because $\S$ is a subsemilattice of $(\F, \cup)$ that is closed downwards in $(\F, \subseteq)$. 
 
 Thus, $\HD(\F)$ contains countably-many principal elements and uncountably-many non-principal elements, 
 $$\text{PHD}(\F) = \{ \S_A: A \subseteq \N \text{ finite}  \} \quad \text{and} \quad \HD(\F) \bs \text{PHD}(\F)= \{ \S_B: B \subseteq  \N \text{ infinite}  \}.$$
 As $\S_B \cap \S_C = \S_{B \cap C}$, $(\HD(\F), \cap ) \cong ({\cal P}(\N), \cap)$. Thus, $ \ds \Delta(A^*_{\F}(G)) = \bigcupdot_{B \subseteq \N} G_{\S_B}$,   with product  $$ G_{\S_B} G_{\S_C} \subseteq G_{\S_{B\cap C}} \quad  \text{ given by }  \quad (s_A)_{A  \in \S_B} (t_{A'})_{A'  \in \S_C} = (s_At_A)_{A  \in \S_{B\cap C}},$$
 is a semilattice of disjoint groups with compatible central identities, graded over $({\cal P}(\N), \cap)$. 
 


\eex

\bex \label{CofiniteExample} \rm Let ${\cal C} = \{ A \subseteq \N: \N \bs A \text{ is finite}\}$. Then  
${\cal C} $ is an example of an infinite upper subsemilattice of $({\cal P}(\N), \cup, \subseteq)$ such that $\HD({\cal C}) = \PHD({\cal C})$: taking $\S$ in $\HD({\cal C})$ and $A_0$ in $\S$, the finite union $A=\cup \{ C \in \S: A_0 \subseteq C\}$ belongs to $\S$ and for any $B$ in $\S$, $A_0 \cup B \in \S$; hence,  $B  \in \S_A$.  Since $\S_A \cap \S_B = \S_{A \cap B}$, the lower semilattices $(\HD({\cal C}), \cap)$ and $({\cal C}, \cap)$ are isomorphic. Thus, $\Delta(A_{\cal C}^*(G)) = \bigcup_{A \in {\cal C}} G_{\S_A}$ is a semilattice of disjoint groups, graded over the lower semilattice $({\cal C}, \cap)$. If, as in Example \ref{Finite Boolean Alg Example}, one identifies $G_{\S_A}$ with $G_A$, the product in  $\Delta(A_{\cal C}^*(G)) = \bigcup_{A \in {\cal C}} G_A$  is described by (\ref{Omnibus product in spectrum of G_A G_B Eqn}). We note  that  $({\cal C}, \subseteq)$ fails to have a minimum element but  $\Delta(A_{\cal C}^*(G))$ is nonetheless  a  semigroup (see Proposition \ref{Delta(A)SemigpProp}.1) and, relevant to Remark \ref{Final Remark}.1, each $A$ in ${\cal C}$ has only finitely many upper bounds.  
\eex 

To further illustrate the scope of our omnibus construction, we give two more examples where  
$A_\E^*(G)$ and $\Delta(A_\E^*(G))$ are graded over infinite countable chains. 

\bex  \rm For $n \in \N$, let $L_n= \{ n, n+1, \ldots\}$. Then ${\cal L} = \{ L_n: n \in \N\}$ is another example of an infinite upper subsemilattice of $({\cal P}(\N), \cup, \subseteq)$ such that $\HD({\cal L}) = \PHD({\cal L})$, ${\cal L}$ has no minimum element, and each $L_n$  has only finitely many upper bounds in ${\cal L}$.  In this case, $({\cal L}, \cup)$ is an infinite chain isomorphic with $(\N, \min)$ and $(\HD({\cal L}), \cap) \cong (\N, \max)$. Letting $\tau_n = \tau_{L_n}$, $$A_{\cal L}^*(G) = \ell^1-\bigoplus_{n \in \N} A_{\tau_n}(G) \quad \text{where}  \quad  A_{\tau_m}(G)  A_{\tau_n}(G) \subseteq  A_{\tau_{m  \wedge n}}(G)$$ and $\Delta(A_{\cal L}^*(G)) = \bigcupdot_{n \in \N} G_{\S_{L_n}} \cong \bigcupdot_{n \in \N} G_{L_n}$ is a semilattice of groups graded over the infinite chain $(\N, \max)$: $G_{L_m} G_{L_n}\subseteq G_{L_{m\vee n}}$. 
\eex

\bex \rm For $n \in \N$, let $K_n= \{ 1, 2, \ldots, n \}$,  ${\cal K} = \{ K_n: n \in \N\}$. Then ${\cal K}$ is an infinite chain in $({\cal P}(\N), \cup, \subseteq)$   isomorphic with $(\N, \max)$.   In this case  $\HD({\cal K}) = \PHD({\cal K}) \cup \{{\cal K}\}$  and $(\HD({\cal K}), \cap) \cong (\N \cup \{ \infty\} , \min)$. Letting $\gamma_n = \tau_{K_n}$, $$A_{\cal K}^*(G) = \ell^1-\bigoplus_{n \in \N} A_{\gamma_n}(G) \quad \text{where}  \quad  A_{\gamma_m}(G)  A_{\gamma_n}(G) \subseteq  A_{\gamma_{m  \vee n}}(G)$$ and, letting $\S_{K_\infty} = {\cal K}$,   $\Delta(A_{\cal K}^*(G)) = \bigcupdot\{G_{\S_{K_n}}: n \in \N \cup\{\infty\}\}$ is a semilattice of groups graded over the infinite chain $(\N \cup \{\infty\} , \min)$.

\eex

Given any finite chain $(C, \leq)$, it will be obvious that one can choose $\E$ so that 
$A_\E^*(G)$ and $\Delta(A_\E^*(G))$ are respectively graded over $(C, \min)$ and $(C, \max)$, or vice versa.

\section{Fusions of compatible maps} 

\noindent Let  $S = \bigcupdot_{z \in Z} G_z$ be a semilattice of disjoint semigroups with compatible central identities $1_z$,  $E_\alpha$ a subset of a group $H$, and 
$ \ds \alpha: E_\alpha \ra S = \bigcupdot_{z \in Z} G_z.$
For each $z \in Z$, let 
\bi 
\item $\fza=\{h \in E_\alpha: \alpha(h) \in G_z\}$;
\item $\kza = \bigcupdot_{w\geq z} \fwa$; 
\item $\bza  : \kza \ra G_z : h\mapsto \alpha(h) 1_z$ whenever $\kza \neq \emptyset$. 
\ei 
We will say that $\alpha$ is a \it fusion \rm of the maps  $\bza : \kza \ra G_z$ $(z \in Z$), which are \it compatible \rm  with $(Z,\leq)$. 

\bigskip 

For the semilattices  of semigroups $S = \bigcupdot_{z \in Z} G_z$ considered herein, and in \cite{Sto}, the lower semilattice $(Z, \cdot, \leq)$ will also be a complete upper semilattice with a maximum element. In this paper, each $G_z$ will be a group, but this is not always the case in \cite{Sto} where the results in this section are also used.  Our choice of terminology was motivated by parts (a) and (b) of Propostition \ref{Fusions of compatible maps Prop}.  Part 2 of the proposition  provides a recipe for producing examples of fusions of compatible maps from a collection of maps $\beta_z:K_z  \ra S$ defined on an admissible family of subsets $K_z$ of $H$. 
 
 \bp \label{Fusions of compatible maps Prop} \rm Let $S = \bigcupdot_{z\in Z} G_z$ be a semilattice of semigroups with compatible, central identities $1_z$ in $G_z$, and let $H$ be a group. 
 \bi \item[1.] Let $  \alpha: E_\alpha \subseteq H \ra S = \bigcupdot_{z \in Z} G_z.$
 Then
 \bi 
\item[(a)] \it Fusion property: \rm $\ds \alpha= \bigcupdot_{z\in Z}\bza\large{\vert}_{\fza} : E_\alpha= \bigcupdot_{z\in Z} \fza \ra S = \bigcupdot_{z\in Z} G_z;$
\item[(b)] \it Compatibility property: \rm for $z \leq w$, $\kza \supseteq \kwa$ and for $h \in \kwa$, $\bwa(h) 1_z= \bza(h).$
\ei
Moreover, if $(Z, \leq) $  is a complete upper semilattice,  then 
\bi \item[(c)] for each nonempty subset $Q$ of $Z$, $\ds \bigcap_{z \in Q} \kza = K_{\vee Q}^\alpha$.
\ei 
 \item[2.] Suppose further that  $(Z, \leq)$ is a complete upper semilattice. For each $z \in Z$,  suppose that $K_z$ is a subset of $H$  and, when $K_z$ is nonempty,   $\beta_z : K_z \ra G_z$  is a mapping such that \bi  \item[(b')] for  $z \leq w$, $K_z \supseteq K_w$ and for $h \in K_w$, $\beta_w(h) 1_z= \beta_z(h);$ and 
 \item[(c')] for each  nonempty subset $Q$ of $Z$, $\ds \bigcap_{z \in Q} K_z = K_{\vee Q}$. 
\ei Putting 
 $$F_z:= K_z \backslash \bigcup_{w>z}K_w \quad \text{(possibly empty),}$$
 $\{F_z: z \in Z\}$ partitions $E_\alpha := \bigcup_{z \in Z} K_z$, and if we define $\alpha$ via \bi \item[(a')] $\ds  \alpha= \bigcupdot_{z\in Z}\beta_z \large{\vert}_{F_z} : E_\alpha= \bigcupdot_{z\in Z} F_z \ra S = \bigcupdot_{z\in Z} G_z, $ \ei
 then for each $z \in Z$, 
 $$F_z = \fza, \quad K_z = \kza \quad \text{and} \quad \beta_z = \bza.$$
(Thus, $\alpha$ is the fusion of the maps $\beta_z:K_z \ra G_z$, which are compatible with $Z$.) 
\ei 
\ep

\begin{proof}  1. It is easy to check that properties (a) and (b) hold. For (c), assume that $(Z, \leq)$ is a complete upper semilattice, and let $Q$  be a nonempty subset of $Z$. From (b), $K_{\vee Q}^\alpha \subseteq \bigcap_{z \in Q} K_z^\alpha$ is clear. Take $h \in \bigcap_{z \in Q} K_z^\alpha$ and let $z_0$ be the unique element of $Z$ such that $h \in F^\alpha_{z_0}$. As $h \in \kza = \bigcupdot_{w \geq z} \fwa$ for each $z \in Q$, we must have $z_0\geq z$ for each $z \in Q$; hence, $z_0 \geq \vee Q$, giving $h \in K_{\vee Q}^\alpha$. 
\smallskip 

\noindent 2.
If $z \vee w > z$ in $Z$, then
$$F_z = K_z \backslash \cup_{z'>z}K_{z'} \subseteq K_z \backslash K_{z \vee w} =  K_z \backslash (K_z \cap K_w) =  K_z \backslash K_w \subseteq K_z \bs F_w,$$
where we have used property (c').  So, if it happens that $F_z \cap F_w$ is nonempty, then $z \vee w = z$ and, by symmetry, $w = z\vee w  =z$. Thus, the sets $\{F_z: z \in Z\}$ are pairwise disjoint and $\bigcupdot_{z \in Z} F_z \subseteq \bigcup_{z \in Z} K_z = E_\alpha$. On the other hand, if  $h \in E_\alpha$ and  we let $Q= \{ z \in Z: h \in K_z\}$, then $h \in \bigcap_{z \in Q} K_z = K_{\vee Q}$ and $h \notin K_w$ for any $w > \vee Q$; so $ h \in F_{\vee Q}$. Thus, $E_\alpha = \bigcupdot_{z \in Z} F_z$ and if we define $\alpha$ via (a'), we immediately obtain $F_z = \fza$ for each $z \in Z$. 

We claim that $\kza = K_z$, where $ \kza = \bigcupdot_{w\geq z} \fwa = \bigcupdot_{w\geq z} F_w$. For $w \geq z$, we know from (b') that $F_w \subseteq K_w \subseteq K_z$, so $\kza \subseteq K_z$. Let $h \in K_z \subseteq E_\alpha = \bigcupdot_{w \in Z} F_w$ and take $w_0 \in Z$ such that $h \in F_{w_0} = K_{w_0} \bs \bigcup_{w > w_0} K_w$. Then $h \in K_{w_0} \cap K_z = K_{w_0 \vee z}$, so we cannot have $w_0 \vee z > w_0$. Thus, $h \in F_{w_0} \subseteq \bigcupdot_{w\geq z} F_w = \kza$ and we have the claim.  

To see that $\beta_z = \bza$, let $h \in K_z = \kza$ and take the unique $w$ in $Z$ with $w \geq z$ such that $h \in \fwa = F_w$. From the definition of $\alpha$, $\alpha (h) = \beta_w(h)$; (b') yields $h \in K_w \subseteq K_z$  and $$\bza(h) = \alpha(h) 1_z = \beta_w(h) 1_z = \beta_z(h).  \qedhere$$

\end{proof} 

\br  \rm  Though it is not needed, observe that the set-containment part of condition (b') follows from (c'):  for $z \leq w$, $K_w = K_{w\vee z} = K_w \cap K_z \subseteq K_z$.
\er 

\bd \label{Def - Fusion of Lambda maps} \rm  If $\Lambda$ is a property of a pair $(K,\beta)$ where $K$ is a subset of a group $H$ and $\beta$ is a mapping of $K$ into a type of semigroup, we will say that $\alpha: E_\alpha \subseteq H \ra S= \bigcupdot_{z \in Z} G_z$ is a fusion of property-$\Lambda$ maps, compatible with $Z$, if $\bza: \kza\ra G_z$ has property $\Lambda$ whenever $\kza \neq \emptyset$.  
\ed 
For example, $\alpha: E_\alpha \subseteq H \ra S= \bigcupdot_{z \in Z} G_z$
is a fusion of group homomorphisms (affine maps), compatible with $Z$, if, whenever $\kza \neq \emptyset$, 
$\kza$ is a subgroup of (coset in) $H$ and 
 $\bza: \kza \ra G_z$ is a homomorphism (affine mapping) into a subgroup of $G_z$.

\br \rm  When $(Z, \leq)$ is a finite chain and each $G_z$ is a group, Definition \ref{Def - Fusion of Lambda maps} agrees with the definitions of compatible fusions of homomorphisms/affine maps found in \cite{Sto21}. 
\er

\bc \label{Fusions of homs - aff maps Cor} \rm   Let $S = \bigcupdot_{z \in Z} G_z$ be a semilattice of semigroups with compatible, central identities. 
\bi \item[1.] If $\alpha: H_\alpha \subseteq H \ra S$ is a fusion of group homomorphisms (affine maps), compatible with $Z$,  then $H_\alpha$ is a subgroup of (coset in) $H$. 
\item[2.]  Suppose that the hypotheses (and conclusions) of Proposition \ref{Fusions of compatible maps Prop}.2 are satisfied. Suppose further that when $K_z$ is nonempty, \bi  \item $K_z$ is a subgroup of (coset in) $H$; and 
\item $\beta_z: K_z \ra G_z$ is a homomorphism (affine map) into a subgroup of $G_z$.
\ei 
Then $$\alpha= \bigcupdot_{z\in Z}\beta_z \large{\vert}_{F_z} : E_\alpha= \bigcupdot_{z\in Z} F_z \ra S = \bigcupdot_{z\in Z} G_z $$
is a fusion of homomorphisms (affine maps), compatible with $Z$. 
\ei 
\ec 

\begin{proof} Part 1: In the case that $\alpha$ is a fusion of affine maps, let $g, h, k \in E_\alpha$, say $g \in K_v^\alpha$, $h \in \kwa$,  $k \in \kza$. By Proposition  \ref{Fusions of compatible maps Prop}.1(b), $g, h, k \in K_v^\alpha \cup \kwa \cup \kza \subseteq K_{vwz}^\alpha$, so $gh^{-1}k \in K_{vwz}^\alpha \subseteq H_\alpha$.  Part 2: By Proposition   \ref{Fusions of compatible maps Prop}.2, $\kza = K_z$ and $\bza =\beta_z$.
\end{proof} 

\bex \label{MainFusionExample} \rm  Let $H$ be an infinite  group and let $G = H^\N$ be given the discrete topology. Let $\F = \{ A \subseteq \N: A \text{ is finite}\}$,  $A_\F^*(G) = \ell^1-\bigoplus_{A \in \F} A_{\tau_A}(G)$ the associated GSA considered in Example \ref{P(N) Example}. As noted there, $$\Delta(A_\F^*(G)) = \bigcupdot_{B \in {\cal P}(\N)} G_{\S_B}$$ 
is a semilattice of groups, graded over the lower semilattice $({\cal P}(\N), \cap, \subseteq)$, with compatible identities $e_{\S_B} = (e_A)_{A \in \S_B}$, 
where 
$\S_B = \{ A \in {\cal P}(\N): A \subseteq B \text{ is finite}\}$ and $G_{\S_B} = \varprojlim_{A \in \S_B} G_A.$

Taking $K$ to be any non-trivial group, we will use Corollary \ref{Fusions of homs - aff maps Cor} to construct a fusion of homomorphisms, compatible with $({\cal P}(\N), \cap, \subseteq)$,
$\alpha:K^\N \ra  \bigcupdot_{B \in {\cal P}(\N)} G_{\S_B}$, that maps non-trivially into each $G_{\S_B}$: 
 For each subset $B$ of $\N$, let 
$$K_B = \{ s= (s_n) \in K^\N: s_n = e_K \text{ for each } n \in B\}.$$
Clearly, each $K_B$ is a subgroup of $K^\N$, $K_B \supseteq K_C$ for $B \subseteq C$, and if ${\cal Q} \subseteq {\cal P}(\N)$, then $\bigvee {\cal Q}  = \bigcup {\cal Q} $ and $\bigcap_{B \in {\cal Q} } K_B = K_{\bigcup {\cal Q} }$. Let $\beta: K^\N \ra G_{\S_\N}$ be any homomorphism (e.g., the trivial one). For $B \in {\cal P}(\N)$, 
$$\beta_B: K_B \ra G_{\S_B}: s \mapsto \beta(s) e_{\S_B}$$ is a homomorphism because $e_{\S_B}$ is a central idempotent in $\Delta(A_\F^*(G))$. If $B \subseteq C$ and $s \in K_C$, then $$\beta_C(s) e_{\S_B} = \beta(s) e_{\S_C}e_{\S_B} = \beta(s) e_{\S_B} = \beta_B(s),$$
so the conditions of Corollary \ref{Fusions of homs - aff maps Cor}.2 are satisfied.  For $B \subseteq \N$, 
$$F_B =  K_B \bs \bigcup_{C \supsetneq B} K_C = \{ s \in K^\N: s_n = e_K \text{ for each } n \in B \text{ and } s_n \neq e_K \text{ for each } n \in \N \bs B \},$$
with $F_\N = K_\N = \{e_{K^\N}\}$ and $F_\emptyset = \{ s \in K^\N: s_n \neq e_K \text{ for each } n \in \N\}$. By Proposition \ref{Fusions of compatible maps Prop}.2 and Corollary \ref{Fusions of homs - aff maps Cor}.2, the nonempty sets $F_B$ partition $K^\N  = \bigcup_{B \subseteq \N} K_B$, and   $$\alpha= \bigcupdot_{B\subseteq \N}\beta_B \large{\vert}_{F_B} : K^\N = \bigcupdot_{B\subseteq \N} F_B \ra \Delta(A_\F^*(G))  = \bigcupdot_{B \subseteq \N} G_{\S_B} $$
is a fusion of group homomorphisms, compatible with $({\cal P}(\N), \cap, \subseteq)$, such that for each $B \in {\cal P}(\N)$,  $F_B^\alpha = F_B$, $K_B^\alpha = K_B$ and $\beta_B^\alpha = \beta_B$. Note that  $\alpha$ maps nontrivially into each $G_{\S_B}$ because each $F_B^\alpha = F_B$ is nonempty. Thus, $\alpha$ is an example of a fusion of compatible homomorphisms mapping nontrivially into infinitely many chains and antichains of subgroups $G_{\S_B}$ for $B$ in $({\cal P}(\N), \cap, \subseteq)$, including uncountably many $G_{\S_B}$ where $\S_B \in \HD(\F)$ is nonprincipal. As we will see in the next section, 
$j_\alpha: A_\F^*(G) \ra B(K^\N): u \mapsto u \circ \alpha$ is a completely positive algebra homomorphism; cf. \cite[Theorem 4.15 and Example 4.16]{Sto21}.

\eex 

\bex \label{FusionP([n])Example}  \rm Let $H$ be an infinite  group and let $G = H^\N$ be given the discrete topology. Consider $\E_n = \{ A \in {\cal P}(\N): A \subseteq [n]\}$ and its  associated GSA, $A_{\E_n}^*(G)$, introduced in  Example \ref{Finite Boolean Alg Example}. As noted there, $\Delta(A_{\E_n}^*(G)) = \bigcupdot_{A \subseteq [n]} G_A$ 
is a semilattice of groups, graded over the lower semilattice $({\cal P}([n]), \cap, \subseteq)$, with compatible central identities $e_A$.  We proceed like in Example \ref{MainFusionExample}: Let $K$ be a nontrivial group, $\beta:K^n \ra G_{[n]}$ a homomorphism,  $K_A = \{ s \in K^n: s_m = e_K \text{ for } m \in A\},$
and $\beta_A: K_A \ra G_A : s \mapsto \beta(s) e_A$. By Corollary \ref{Fusions of homs - aff maps Cor}, 
$$\alpha= \bigcupdot_{A\subseteq [n]}\beta_A \large{\vert}_{F_A} : K^n = \bigcupdot_{A\subseteq [n]} F_A \ra \Delta(A_{\E_n}^*(G))  = \bigcupdot_{A \subseteq [n]} G_A $$ 
is a fusion of  homomorphisms, compatible with $({\cal P}([n]), \cap, \subseteq)$, such that each $K_A^\alpha =K_A$, $\beta_A^\alpha = \beta_A$. Observe that each $F_A^\alpha = F_A = \{ s \in K^n: s_m = e_K \text{ for } m \in A \text{ and } s_m \neq e_K \text{ for  } m \in [n] \bs A \}$ is non-empty, so  $\alpha$ maps nontrivially into  $G_A$ for each $A \subseteq [n]$.  Figure 2 provides  an illustration of this fusion of compatible homomorphisms when $[n]=[2]= \{1,2\}$.  \eex

\begin{figure} \label{Figure2}
\centering
\begin{tikzpicture}[
  every node/.style={font=\small},
  region/.style={draw, thick},
  squareObj/.style={draw, thick, minimum size=1.1cm},
  map/.style={->, thick},
  x=1cm, y=1cm
]


\draw[region] (0,0) circle (2.4);

\draw[region] (-0.8,0) circle (1.4);
\draw[region] ( 0.8,0) circle (1.4);

\node (K2) at (-1.4,0) {$\scriptscriptstyle{K_{\{2\}}}$};    
\node (K1) at ( 1.4,0) {$\scriptscriptstyle{K_{\{1\}}}$};    
\node (K0) at ( 0,0) {$\scriptscriptstyle{K_{\{1,2\}}}$};      
\node (H) at ( 0,-1.5) {$\scriptscriptstyle{K_\emptyset}$};   


\node[squareObj] (Q1) at (6, 2) {\ \ $\scriptscriptstyle{G_{\{1,2\}} \ \ }$};
\node[squareObj] (Q2) at (4.6,0) {$\scriptscriptstyle{G_{\{1\}}}$};
\node[squareObj] (Q3) at (7.4,0) {$\scriptscriptstyle{G_{\{2\}}}$};
\node[squareObj] (Q4) at (6,-2) {$\scriptscriptstyle{G_\emptyset=G^{ap}}$};


\draw[map]
  (K2)
  .. controls +(300:5) and +(290:5) ..
  node[midway, above left] {$\scriptscriptstyle\beta_{\{2\}}$}
  (Q3.south east);
\draw[map] (K1) to[bend left=15] node[midway, above right] {$\scriptscriptstyle\beta_{\{1\}}$} (Q2.west);
\draw[map] (0,1) to[bend left=10]  node[midway, above]       {$\scriptscriptstyle\beta_{\{1,2\}}$} (Q1.west);
\draw[map] (H) to[bend right=10] node[midway, above]       {$\scriptscriptstyle\beta_{\emptyset}$} (Q4.west);
\draw[map] ([xshift=20pt]Q1.south west) to node[midway, left] {$\scriptscriptstyle (\cdot)  e_{\{1\}}$} (Q2.north);
\draw[map] ([xshift=-20pt]Q1.south east) to node[midway, right] {$\scriptscriptstyle(\cdot)  e_{\{2\}}$} (Q3.north);
\draw[map] (Q2.south) to node[midway, left] {$\scriptscriptstyle (\cdot)  e_{\emptyset}$} ([xshift=10pt]Q4.north west);
\draw[map] (Q3.south) to node[midway, right] {$\scriptscriptstyle (\cdot)  e_{\emptyset}$} ([xshift=-10pt]Q4.north east);

\end{tikzpicture}
\vspace{-.7cm}
\caption{ $\ds \alpha= \bigcupdot_{A\subseteq \{1,2\}}\beta_A \large{\vert}_{F_A} : K^2 = \bigcupdot_{A\subseteq \{1,2\}} F_A \ra \Delta(A_{\E_2}^*(G))  = \bigcupdot_{A \subseteq \{1,2\}} G_A$ \\ from Ex. \ref{FusionP([n])Example} when $n=2$: a fusion of homomorphisms, compatible with $({\cal P}(\{1,2\}), \cap)$}
\end{figure}

\section{Homomorphisms of spine algebras}

For any GSA  $A_\B^*(G)$,  we will characterize all completely positive and completely contractive homomorphisms of $A_\B^*(G)$ into $B(H)$; when $G$ is amenable, we will characterize every such   completely bounded homomorphism. Our results are new, even for the full spine algebra $A^*(G)$ including, in the  (completely) positive and contractive cases,   when $G$ is abelian. When $G$ and $H$ are abelian,  Inoue provided a description of the (automatically completely bounded) homomorphisms of $ A^*(G)$ into   $B(H)$ \cite[Theorem 4.1]{Inoue} that is related to statement 2 of our Theorem \ref{CPCCCBHomomThm} below. 

Let $\vp: A \ra B(H)$ be a homomorphism, where $A= A_\pi$ is a closed translation-invariant subalgebra of $B(G)$.  Since $B(H)$ is semisimple, $\vp$ is  bounded and, by \cite[Proposition 1.8]{Sto21}, there is an open subset $E$ of $H$ and a  continuous map $\alpha: E \ra \Delta(A)$ such that $\vp = j_\alpha$: for $u\in A$ and $h \in H$, 
$$\displaystyle{\vp(u)(h) = j_\alpha(u)(h) = \left \{ \begin{array}{ll}
                          \l \alpha(h), u \r  & \mbox{if  $ h \in  E$}\\
                          \ \ \  0 & \mbox{if  $ h \in H \bs E$.}
                            \end{array}
                 \right. }   $$ 
 Here, $E = \{ h \in H : \vp^*(h) \neq 0\}$ where $\vp^*: W^*(H) \ra VN_\pi$  is the dual operator and $\alpha$ is the restriction of $\vp^*$ to $E$. Equivalently, there is a continuous map $\alpha_0: H  \ra \Delta(A)\cup \{0\}$ such that $E = \{ h \in H: \alpha (h) \neq 0 \}$ is open and $\vp(u) (h) = j_{\alpha_0}(u)(h) = \l \alpha_0(h), u \r $ for $ u \in A$, $h \in H$. The set $E$ and the maps $\alpha$, $\alpha_0$ are uniquely determined by $\vp$.

\blem  \label{Homom into ellinfinity Lemma} \rm   Let $A$ be a closed translation-invariant  subalgebra of $B(G)$, $\gamma:E\subseteq H \ra \Delta(A)$ and define $\gamma_0: H \ra \Delta(A) \cup \{ 0\}$ by $\gamma_0 \vert_E = \gamma$, $\gamma_0 \vert_{H\bs E} = 0$.  Then  $j_\gamma = j_{\gamma_0}: A \ra \ell^\infty(H)$ is a contractive homomorphism. If $j_\gamma$ maps into $B(H)$, then $j_\gamma: A \ra B(H)$ is a bounded homomorphism, $E$ is an open subset of $H$ and $\gamma$ is continuous.
\elem  

\begin{proof}  The first statement is easy to check.  If $j_\gamma$ maps into $B(H)$, then $j_\gamma: A \ra B(H)$ is an automatically bounded homomorphism, and $E$, $\gamma$ are the open set and continuous map uniquely determined by $j_\gamma$, identified in  \cite[Proposition 1.8]{Sto21} (and described above).
\end{proof}

Throughout the remainder of this section, $\A = A_\B^*(G)$ is GSA, where $\B = \B_\A$,  and  $G_\B^* = \Delta(\A) = \bigcupdot_{\S \in \HD(\B)} G_\S$.    Let $\alpha: E_\alpha \subseteq H \ra G_\B^*$. As before, and throughout the remainder of this section, we will use the notation: \bi \item  $F_\S = F_\S^\alpha = \{ h \in E_\alpha: \alpha (h) \in G_\S\}$;  \item  $K_\S = K_\S^\alpha = \bigcupdot_{\S' \supseteq \S} F_{\S'}$; 
\item $\beta_\S = \beta_\S^\alpha: K_\S \ra G_\S: h \mapsto \alpha(h) e_\S$ when $K_\S \neq \emptyset$. 
\ei 
If $\tau\in\B$ and  $K_{\S_\tau} \neq \emptyset$, let $$\beta_\tau = \widetilde{\beta}_{\S_\tau}: K_{\S_\tau} \ra G_\tau: h \mapsto \beta_{\S_\tau}(h)^\sim,$$
where $G_{\S_\tau} \ra G_\tau: s =(s_\sigma)_{\sigma \in \S_\tau} \mapsto \widetilde{s} := s_\tau,$ a topological group isomorphism.  Let $$\S_\alpha  =\{ \tau \in \B: K_{\S_\tau} \neq \emptyset\}.$$ 

\br  \rm  At this point, we are not assuming that $\S \cap \S' \neq \emptyset$ for $\S, \S' \in \HD(\B)$, i.e., we are not assuming that $G_\B^* = \Delta(A_\B^*(G))$ is a semigroup. \er

 \bp  \label{FirstHMProp}  \rm Let $\alpha : E_\alpha \subseteq H \ra G_\B^* = \bigcupdot_{\S \in \HD(\B)} G_\S$.  
 
 \begin{enumerate}  \item The map $j_\alpha$ is a contractive homomorphism of $A_\B^*(G)$ into $\ell^\infty (H)$ such that $j_\alpha\large{\vert}_{A_\tau(G)} = j_{\beta_\tau}$ for $\tau \in \S_\alpha$ and $j_\alpha\large{\vert}_{A_\tau(G)} =0 $ for $\tau \in \B \bs \S_\alpha$. Thus, 
 \beq \label{jalpha vs jbetatau eqn}  j_\alpha u = \sum_{\tau \in \S_\alpha}   j_{\beta_\tau} u_\tau \quad \text{for } u = \sum_{\tau \in \B} u_\tau \in A_\B^*(G) =  \ell^1-\bigoplus_{\tau \in \B} A_\tau(G).  \eeq 
 \item The following statements are equivalent:  \bi \item[(i)] $j_\alpha$ is a (automatically bounded, homomorphic) mapping 
 of $A_\B^*(G)$ into $B(H)$;
 \item[(ii)]  for each $\tau$ in $\S_\alpha$, $j_{\beta_\tau}$ is a (automatically bounded, homomorphic) mapping 
 of $A_\tau(G)$ into $B(H)$ with $\sup_{\tau \in \S_\alpha} \| j_{\beta_\tau}\| < \infty$.
 \ei When these equivalent statements hold, (\ref{jalpha vs jbetatau eqn}) holds with convergence in $B(H)$, $\|j_\alpha\| = \sup_{\tau \in \S_\alpha} \|j_{\beta_\tau}\|$,  $E_\alpha$ is open, $\alpha$ is continuous, and for each $\tau$ in $\S_\alpha$, $K_{\S_\tau}$ is open and $\beta_\tau$ is continuous. 
 \item For $\S$ in $\HD(\B)$, $K_\S = \bigcap_{\tau \in \S} K_{\S_\tau}$ and $$\beta_\S: K_\S \ra G_\S = \varprojlim_{\tau \in \S} G_\tau: h \mapsto (\beta_\tau(h))_{\tau \in \S} \quad \text{when } K_\S \neq \emptyset.$$
 \end{enumerate} 
 \ep 
 
 \begin{proof}  1. By Lemma \ref{Homom into ellinfinity Lemma}, $j_\alpha: \ABSG \ra \ell^\infty(H)$  is a contractive homomorphism. Let $h \in E_\alpha$, say $h \in F_\S$, so $\alpha(h) = (\alpha(h)_\sigma)_{\sigma \in \S} \in G_\S = \varprojlim_{\sigma \in \S} G_\sigma$. For $\tau$ in $\S$, $h \in K_{\S_\tau}$ and $\beta_{\S_\tau}(h) = \alpha(h)e_{\S_\tau}= (\alpha(h)_\sigma)_{\sigma \in \S_\tau}$. Hence, \beq \label{FirstHMPropEqn2} \beta_\tau(h) = \widetilde{\beta}_{\S_\tau}(h) = \alpha(h)_\tau \quad \text{for } h \in F_\S \text{ and }  \tau \in \S. 
 \eeq
 Let $\tau \in \B$, $u_\tau \in A_\tau(G)$, $h \in E_\alpha$ with $h \in F_\S$, as above.  Writing $\alpha(h) = (\alpha(h)_\sigma)_{\sigma \in \B}$, where $\alpha(h)_\sigma = 0$ for $\sigma \in \B \bs \S$, 
 \beq \label{FirstHMPropEqn1}    j_\alpha u_\tau(h) = \l \alpha(h), u_\tau\r= \l \alpha(h)_\tau, u_\tau\r. 
 \eeq  
 If $h \in E_\alpha \bs K_{\S_\tau}$, then $\tau \notin \S$, so $j_\alpha u_\tau(h) = 0$; hence, $j_\alpha u_\tau\vert_{H\bs K_{\S_\tau}}  = 0$. Thus,  $j_\alpha\vert_{A_\tau(G)} = 0$ for $\tau \notin \S_\alpha$. Suppose that $\tau \in \S_\alpha$. For $h \in H \bs K_{\S_\tau}$, we have observed that $j_\alpha u_\tau (h) = 0 = j_{\beta_\tau}u_\tau (h)$, so take $h \in K_{\S_\tau}$. If $h \in F_\S$, then $\tau \in \S_\tau \subseteq \S$, so 
 $$  j_\alpha u_\tau(h) =  \l \alpha(h)_\tau, u_\tau\r = \l \beta_\tau(h), u_\tau\r = j_{\beta_\tau}u_\tau(h)$$ 
 by (\ref{FirstHMPropEqn2}) and (\ref{FirstHMPropEqn1}). 
 We conclude that $j_\alpha\vert_{A_\tau(G)} = j_{\beta_\tau}$ for $\tau \in \S_\alpha$; (\ref{jalpha vs jbetatau eqn}) follows.  
 
 Statement 2 follows readily from statement 1 and Lemma \ref{Homom into ellinfinity Lemma}.  Let $\S \in \HD(\B)$. Let $h \in K_\S$, say $h \in F_{\S'}$. Then $\S \subseteq \S'$, so  $S_\tau \subseteq \S'$ for each $\tau \in \S$, whence $h \in K_{\S_\tau}$. Thus, $\K_\S \subseteq \bigcap_{\tau \in \S}K_{\S_\tau}$. On the other hand, take $h \in \bigcap_{\tau \in \S}K_{\S_\tau} $.  Again, $h \in F_{\S'}$ for a unique $\S'$ in $\HD(\B)$ such that $\S' \supseteq \S_\tau$ for each $\tau \in \S$, so  $\S \subseteq  \S'$. We conclude that $h \in K_\S$.  
 
 Finally, take $h \in K_\S$, say $h \in F_{\S'}$ where $\S'\supseteq \S$. By (\ref{FirstHMPropEqn2}),  $\alpha(h)_\tau = \beta_\tau(h)$ for each $\tau \in \S'$, so $\beta_\S(h) = \alpha(h)e_\S = (\alpha(h)_\tau)_{\tau \in \S} = (\beta_\tau(h))_{\tau \in \S}$. 
  \end{proof} 
  
 The following lemma summarizes, and adapts to our context, the most important results concerning completely positive, completely contractive and completely bounded homomorphisms from \cite{Daws,Il, Il-Sp, Pham, Sto21}.

  \blem \label{CPCCCBHom on Atau(G) Lemma} \rm  Let $\tau \in \TG$ and let  $\vp_\tau = j_{\beta_\tau}: A_\tau(G) \ra B(H)$ be a nonzero homomorphism, where $\beta_\tau : E \subseteq H \ra \Delta(A_\tau(G)) = G_\tau$. 
  \bi \item[(a)] Then  $\vp_\tau$ is a completely positive (completely contractive) homomorphism if and only if $E$ is an open subgroup of (coset in) $H$ and $\beta_\tau$ is a continuous homomorphism (affine map). 
  \item[(b)] If $E_\tau \in \Omega(H)$ and $\beta_\tau$ is a continuous piecewise affine map, then $\vp_\tau$ is a completely bounded homomorphism. The converse holds when $G_\tau$ is amenable. 
  \ei
  \elem

  \begin{proof}  If $E$ is an open subgroup, coset, member of $\Omega(H)$, and $\beta_\tau$ is respectively a continuous homomorphism, affine mapping, piecewise affine mapping, then $\vp_\tau = j_{\beta_\tau}$ is, respectively, a completely positive, completely contractive, completely bounded homomorphism by Corollary 2.6 and Propositions 2.12, 2.13 of \cite{Sto21}. Conversely, suppose that $\vp_\tau = j_{\beta_\tau}$ is a completely positive/completely contractive/completely bounded homomorphism. Since $j_{\eta_\tau}$ is a completely positive and completely isometric isomorphism, $\vp_\tau \circ j_{\eta_\tau}: A(G_\tau) \ra B(H)$ is also a  completely positive/completely contractive/completely bounded homomorphism and, if we now identify   $G_\tau$ with $\Delta(A(G_\tau))$ and view $\beta_\tau$ as $\beta_\tau: E \ra G_\tau =  \Delta(A(G_\tau))$, it follows from (\ref{Delta(Atau(G)Identitification}) that $\vp_\tau \circ j_{\eta_\tau} = j_{\beta_\tau}$ as well.  When  $\vp_\tau$ is completely positive (resp. completely contractive),  $E$ is therefore an open subgroup (resp. coset) in $H$ and $\beta_\tau$ is a continuous homomorphism (affine mapping) by Proposition 1.2 (resp. Remark 1.3.2) of \cite{Sto21}; alternatively, one can directly employ Corollary 4.3 (resp. Corollary 4.20) of  \cite{Sto21}, or Proposition 5.8 of \cite{Pham}.  If $G_\tau$ is amenable and $\vp_\tau$ is completely bounded, so too is  $\vp_\tau \circ j_{\eta_\tau} = j_{\beta_\tau}$, so $E \in \Omega(H)$ and $\beta_\tau$ is continuous and piecewise affine by  the main result in \cite{Daws}. 
  \end{proof} 
  
  The following statement  will be known, but is included because we do not have a reference. 
   
  \blem \label{CPCCCB Lemma} \rm Let   $T_i: {\cal A} \ra {\cal B}_i$  be a uniformly bounded collection of continuous linear maps between C$^*$-algebras ${\cal A}$ and ${\cal B}_i$, $i \in I$. If each $T_i$ is completely positive/completely contractive/completely bounded with $\sup\|T_i\|_{cb} < \infty$, then $$T : {\cal A} \ra  \B:= \ell^\infty-\bigoplus_{i \in I} {\cal B}_i: x \mapsto (T_i(x))_{i \in I}$$ is  completely positive/completely contractive/completely bounded. 
  \elem 
  
  \begin{proof}  If each $T_i$ is completely positive, then so is $T$, for example by \cite[Corollary IV.3.4]{Tak}: for $x_1, \ldots, x_n$ in $\cal A$ and $y_1, \ldots, y_n$ in $\cal B$, say $y_l = (y_{l,i})_{i\in I}$ for $l =1, \ldots, n$, 
  $$\sum_{k,l =1}^n y_k^* T(x_k^*x_l)y_l = \left(  \sum_{k,l =1}^n y_{k,i}^* T(x_k^*x_l)y_{l,i}    \right)_{i \in I}  \geq 0.$$ 
  
 Suppose that each $T_i$ is completely bounded with $L= \sup\|T_i\|_{cb} < \infty$. Viewing $\B_i$ as a C$^*$-subalgebra of ${\cal L}(\H_i)$, the space of bounded linear operators on the Hilbert space $\H_i$, $\B$ inherits its operator space structure through the containments 
 $$\B = \ell^\infty-\bigoplus_{i \in I} \B_i \subseteq \ell^\infty-\bigoplus_{i \in I} {\cal L}(\H_i) \subseteq {\cal L} (\ov{\H}),$$
  where $\ov{\H} = \bigoplus_{i\in I} \H_i$. Letting $n$ be a positive integer, we will show that $\|T_n\| \leq L$, where $T_n$ is the $n^{\text{th}}$-amplification of $T$, $T_n : M_n({\cal A}) \ra M_n (\B) \subseteq {\cal L}(\ov{\H}^n)$. Let $x= [x_{jk}] \in M_n({\cal A}) $, $\ov{\xi} = (\xi_k)_{k=1}^n \in \ov{\H}^n$, say $\xi_k = (\xi_{k,i})_{i \in I} \in \ov{\H}$. For each $i \in I$, let $\xi_i: = (\xi_{k,i})_{k=1}^n \in \H_i^n$. Noting that $(T_i)_n: M_n({\cal A}) \ra M_n({\cal L}(\H_i))= {\cal L}(\H_i^n)$, we obtain 
  \beqs  \| (T_n x)\ov{\xi}\|^2 & = & \left\| \left( \sum_{k=1}^n T x_{jk}\xi_k\right)_{j=1}^n\right\|^2 = \left\| \left( \left(\sum_{k=1}^n T_i x_{jk}\xi_{k,i} \right)_{i \in I} \right)_{j=1}^n\right\|^2 \\
  & = & \sum_{j=1}^n \sum_{i \in I}   \left\|  \sum_{k=1}^n T_i x_{jk}\xi_{k,i} \right\|^2 =   \sum_{i \in I} \sum_{j=1}^n   \left\|  \sum_{k=1}^n T_i x_{jk}\xi_{k,i} \right\|^2 \\
  & = & \sum_{i \in I}   \left\| \left( \sum_{k=1}^n T_i x_{jk}\xi_{k,i} \right)_{j=1}^n \right\|^2 
  =   \sum_{i \in I}   \left\| (T_i)_n  x (\xi_i)  \right\|^2 \\ 
  & \leq &  \sum_{i \in I}   \| (T_i)_n\|^2  \|x\|^2\| \xi_i  \|^2 \leq L^2 \|x\|^2 \sum_{i \in I} \| \xi_i\|^2 = L^2 \| x\|^2 \| \ov{\xi}\|^2.
  \eeqs
  Hence, $\|T_n \| \leq L$, as needed. 
  \end{proof}

 \bt   \label{CPCCCBHomomThm}  \rm Let $A_\B^*(G)$ be a GSA.
 \bi \item[1.]  The following statements are equivalent: 
 \bi \item[(a)] $\vp: A_\B^*(G) \ra B(H)$ is a nonzero completely positive (completely contractive) homomorphism; 
 \item[(b)] there is an open subset $E_0$ of  $H$ and a continuous map
 $$\alpha: E_0 \ra G_\B^*  = \bigcupdot_{\S \in \HD(\B)} G_\S$$ such that for each $\tau \in \S_\alpha$, $K_{\S_\tau}$ is an open subgroup of (coset in) $H$, $\beta_\tau:K_{\S_\tau} \ra G_\tau$ is a continuous homomorphism (affine map), and $\vp(u) = j_\alpha u$ for each $u$ in $A_\B^*(G)$.
 \ei  
 \item[2.] When $G$ is  amenable as a  locally compact group, the following statements are equivalent: 
  \bi 
   \item[(a)] $\vp: A_\B^*(G) \ra B(H)$ is a nonzero completely bounded homomorphism; 
  \item[(b)] there is an open subset $E_0$ of  $H$ and a continuous map
 $$\alpha: E_0 \ra G_\B^*  = \bigcupdot_{\S \in \HD(\B)} G_\S$$ such that for each $\tau \in \S_\alpha$, $K_{\S_\tau} \in \Omega(H)$, $\beta_\tau:K_{\S_\tau} \ra G_\tau$ is a continuous piecewise affine map with  $\sup_{\tau \in \S_\alpha} \|j_{\beta_\tau}\|_{cb} < \infty$, and $\vp(u) = j_\alpha u$ for each $u$ in $A_\B^*(G)$. \ei  
 \ei 
  When the equivalent conditions in statements 1 and 2 hold,   \beq \label{jalpha vs jbetatau eqnMainCPThm}  j_\alpha u = \sum_{\tau \in \S_\alpha}   j_{\beta_\tau} u_\tau \quad \text{for } u = \sum_{\tau \in \B} u_\tau \in A_\B^*(G) =  \ell^1-\bigoplus_{\tau \in \B} A_\tau(G)  \eeq 
 holds with convergence in $B(H)$. For $\S$ in $\HD(\B)$, \beq  \label{KSbetaS General formula Eqn} K_\S = \bigcap_{\tau \in \S} K_{\S_\tau}, \qquad  \beta_\S: K_\S \ra G_\S = \varprojlim_{\tau \in \S} G_\tau: h \mapsto (\beta_\tau(h))_{\tau \in \S} \eeq  when  $K_\S $ is nonempty, and $\alpha= \bigcupdot_{\S \in \HD(\B)} \beta_\S\vert_{F_\S}$. 
 \et  
 
 \br  \rm  In Theorem \ref{CPCCCBHomomThm}.2, (b) implies (a) holds without assuming that $G$ is amenable. \er 
 
 \begin{proof}    Take $\alpha$ as described in statement 1(b). For $\tau$ in $\S_\alpha$, $j_{\beta_\tau}: A_\tau(G) \ra B(H)$  is a completely positive (completely contractive) homomorphism by Lemma \ref{CPCCCBHom on Atau(G) Lemma}. In both cases, $\sup_{\tau \in \S_\alpha} \| j_{\beta_\tau}\| < \infty$, so by Proposition \ref{FirstHMProp}.2, $j_\alpha:A_\B^*(G) \ra B(H)$ is a homomorphism satisfying  (\ref{jalpha vs jbetatau eqnMainCPThm}), with convergence in $B(H)$.   Letting $\iota_\tau: A_\tau(G) \hookrightarrow A_\B^*(G)$, the dual map $p_\tau=\iota_\tau^*$ is the projection map $VN_{\lambda_\B} = \ell^\infty-\bigoplus_{\sigma \in \B} VN_\sigma \ra VN_\tau$. Since $j_\alpha \circ \iota_\tau = j_{\beta_\tau}$ for $\tau \in \S_\alpha$ and $j_\alpha \circ \iota_\tau =0$ otherwise, 
 $$j_\alpha^* : W^*(H) \ra VN_{\lambda_\B} = \ell^\infty-\bigoplus_{\tau \in \B} VN_{\tau}: x \mapsto (p_\tau(j_\alpha^*(x)))_{\tau \in \B} = (j_{\beta_\tau}^*(x))_{\tau \in \S_\alpha} \oplus (0)_{\tau \in \B \bs \S_\alpha},$$
 which is completely positive by Lemma \ref{CPCCCB Lemma}.  If condition 2(b) is satisfied, the same argument shows that $j_\alpha$ is a  completely bounded homomorphism satisfying (\ref{jalpha vs jbetatau eqnMainCPThm}).    
 
 Conversely, suppose  that  statement 1(a) or 2(a) holds.   Take $E_0$ to be the open subset of $H$ and $\alpha: H_0 \ra G_\B^*$ the  continuous map, uniquely determined by $\vp$, such that $\vp = j_\alpha$ \cite[Proposition 1.8]{Sto21}.     Let $\tau \in \S_\alpha$.  By Proposition \ref{FirstHMProp}.1, $j_{\beta_\tau} = j_\alpha\vert_{A_\tau(G)} = j_\alpha \circ \iota_\tau$, which, since $\iota_\tau$ is completely positive and completely contractive, is completely positive (completely contractive) when statement 1(b) holds, and completely bounded with $\|j_{\beta_\tau}\|_{cb} < \|j_\alpha\|_{cb}$ when statement 2(a) holds. Assuming 1(a), $K_{\S_\tau}$ is an open subgroup of (coset in) $H$ and $\beta_\tau: K_{\S_\tau} \ra G_\tau$ is a continuous homomorphism (affine map) by Lemma \ref{CPCCCBHom on Atau(G) Lemma}. If $G$ is amenable, each $G_\tau$ is also amenable because $\eta_\tau: G \ra G_\tau$ is a continuous dense-range homomorphism. When $G$ is amenable and statement 2(a) holds, Lemma \ref{CPCCCBHom on Atau(G) Lemma} thus gives statement (b). The final statement is contained in Propositions \ref{FirstHMProp}.3 and \ref{Fusions of compatible maps Prop}.1.  \end{proof}

 Recall that  Proposition \ref{Delta(A)SemigpProp}  describes when $G_\B^*$ is  a semigroup, e.g., when $\B$ has a minimum element. In particular, the spectrum $G^*$ of the full spine algebra $A^*(G)$ is always semigroup.

 \bt \label{MainCPCCHomomThm} \rm Suppose that $G_\B^* = \Delta(A_\B^*(G))$ is a semigroup. The following statements are equivalent: 
 \bi \item[(a)] $\vp: A_\B^*(G) \ra B(H)$ is a nonzero completely positive (completely contractive) homomorphism; 
 \item[(b)] there is an open subgroup (coset) $H_0$ in $H$ and a continuous fusion of homomorphisms (affine maps)
 $$\alpha: H_0 \ra G_\B^*  = \bigcupdot_{\S \in \HD(\B)} G_\S,$$ compatible with the lower semilattice $(\HD(\B), \cap, \subseteq)$, such that $\vp(u) = j_\alpha u $ for each $u$ in $A_\B^*(G)$.
 \ei  
   
 
 \et
 
 \begin{proof}    Take $\alpha$ as described in statement (b). For $\S$ in $\HD(\B)$,  $\beta_\S: K_\S \ra G_\S: h \mapsto \alpha(h)e_\S$ is then a continuous homomorphism (affine map) on the subgroup (coset)  $K_\S$, when $K_\S$ is nonempty.  For $\tau$ in $\S_\alpha$, $\beta_\tau=\widetilde{\beta}_{\S_\tau}: K_{\S_\tau} \ra G_\tau$ is hence a homomorphism (affine map),   so  $j_{\beta_\tau}$ maps $A_\tau(G)$ into $B(H_d)$, where $H_d$ is $H$ with the discrete topology, by Lemma \ref{CPCCCBHom on Atau(G) Lemma}.  But $H_0$ is open (and closed) in $H$, and $\alpha$ is continuous, so $j_\alpha u$ is continuous on $H$ for each $u$ in $A_\B^*(G)$. For $\tau \in \S_\alpha$, $j_{\beta_\tau} = j_\alpha\vert_{A_\tau(G)} $ hence maps $A_\tau(G)$ into $B(H_d) \cap CB(H) = B(H)$ by Proposition \ref{FirstHMProp}.1 and \cite[Corollary 2.22]{Kan-Lau}.  By Lemma \ref{Homom into ellinfinity Lemma}, $K_{\S_\tau}$ is thus  open in $H$ and $\beta_\tau$  is  continuous  on $K_{\S_\tau}$ for $\tau \in \S_\alpha$, so   $\vp = j_\alpha$ is a completely positive (completely contractive) homomorphism by Theorem \ref{CPCCCBHomomThm}.

 Suppose  that statement (a) holds, so  $\vp=j_\alpha$ where $\alpha:E_0 \ra G_\B^*$ is a continuous map satisfying the conditions of  Theorem \ref{CPCCCBHomomThm}.1(b).    When $K_\S$ is nonempty, we must show that $K_\S$ is a subgroup of (coset  in) $H$ and $\beta_\S$ is a homomorphism (affine map).   By Theorem  \ref{CPCCCBHomomThm},    $K_{\S_\tau}$ is an open subgroup of (coset in) $H$ and $\beta_\tau:K_{\S_\tau} \ra G_\tau$ is a continuous homomorphism (affine map) for each $\tau \in \S_\alpha$. For $\S$ in  $\HD(\B)$, $K_\S$ and $\beta_\S: K_\S \ra G_\S$ are described by (\ref{KSbetaS General formula Eqn}). Recalling that $K$ is a coset in $H$ exactly when $gh^{-1}k \in K$ for $g, h, k \in K$ and $\beta$ is an  affine mapping on a coset $K$  exactly when $\beta(gh^{-1}k)= \beta(g)\beta(h)^{-1}\beta(k)$ for $g,h,k \in K$, it follows that  $K_\S = \bigcap_{\tau \in \S}K_{\S_\tau}$ is a subgroup (coset)   and  $\beta_\S(h) = (\beta_\tau(h))_{\tau \in \S}$ is a homomorphism (affine map) when $K_\S$ is nonempty.   By Corollary \ref{Fusions of homs - aff maps Cor}, $H_0$ is a subgroup of (coset in) $H$. 
  \end{proof}

\br \rm For non-principal $\S$ in $\HD(\B)$, the above proof shows that the homomorphism (affine map)  $\beta_\S$ is continuous but does not show that the subgroup (coset) $K_\S$ is open. \er

It may happen that $G_\B^* =  \Delta(A_\B^*(G))$ is not a semigroup. (For example, take $\E$ to be a subsemilattice of $({\cal P}(\N), \cup)$  that contains disjoint sets $A$ and $B$, but does not contain the empty set. Then $\S_A \cap \S_B = \emptyset$, so taking $\B = \B_\E$ in the omnibus construction from Section 4, $G_\B^* =\bigcup_{\S \in \HD(\E)} G_\S$ is not a semigroup.) In this case, we can replace $G_\B^*$ with the semigroup $G_\B^*\cup \{0\}$ --- see Proposition \ref{Delta(A)SemigpProp} --- and obtain the following statement using simple modifications to the proof of Theorem \ref{MainCPCCHomomThm}.

\bt  \label{CPCCHomomThmsNoSemigp} \rm Let $A_\B^*(G)$ be any GSA with spectrum $G_\B^* = \Delta(A_\B^*(G))$. The following statements are equivalent: 
 \bi \item[(a)] $\vp: A_\B^*(G) \ra B(H)$ is a nonzero completely positive (completely contractive) homomorphism; 
 \item[(b)] there is  a continuous fusion of homomorphisms (affine maps)
 $$\alpha: H \ra G_\B^* \cup \{0\},$$ compatible with the lower semilattice $(\HD(\B) \cup \{\emptyset\}, \cap, \subseteq)$, such that $\vp(u) = j_\alpha u $ for each $u$ in $A_\B^*(G)$.
 \ei  
\et 
 
\br \rm  \label{Scope of alpha for CP Homs Remark} 1. Let $\vp: A_\B^*(G) \ra B(H)$ be a completely positive homomorphism, where $\Delta(A_\B^*(G)) = G_\B^*$ is a semigroup.  By Theorem \ref{MainCPCCHomomThm}, there is a unique open subgroup $H_0$ and a unique continuous fusion of homomorphisms $\alpha: H_0 \ra G_\B^*$, compatible with the lower semilattice $(\HD(\B), \cap, \subseteq)$, such that $\vp = j_\alpha$. Thus, 
$\vp = \bigcupdot_{\S \in \HD(\B)} \beta_\S\vert_{F_\S}: H_0 = \bigcupdot_{\S \in \HD(\B)} F_\S \ra G_\B^* = \bigcupdot_{\S \in \HD(\B)} G_\S,$
where $\beta_\S:K_\S \leq H_0 \ra G_\S$ is a homomorphism and $F_\S= \{ h \in H_0: \alpha (h) \in G_\S\}$. Theorem \ref{CPCCCBHomomThm} might appear to suggest that $\alpha$ can only map nontrivially into $G_\S$, equivalently $F_\S$ is nonempty, when $\S$  is a  principal element $\S = \S_\tau$ of $\HD(\B)$. In the example below, we observe that there exist completely positive homomorphisms $\vp = j_\alpha$ on $A_\B^*(G)$ where $\alpha$ is a fusion of homomorphisms mapping nontrivially into $G_\S$ for uncountably many nonprincipal $\S$ in $\HD(\B)$. 

\smallskip  

\noindent 2. In \cite{Sto21}, the author considered completely positive homomorphisms $\vp: A \ra B(H)$, where $A$ is a closed translation-invariant subalgebra of $B(G)$ and $\Delta(A)$ is a semilattice of disjoint groups $\bigcupdot_{z\in Z} G_z$ with $(Z, \leq)$  a finite chain. In this scenario, Theorem 4.15 of \cite{Sto21} shows that  $\vp: A \ra B(H)$ is a completely positive homomorphism exactly when $\vp = j_\alpha$ where $\alpha: H_0  \ra \Delta(A) = \bigcupdot_{z\in Z} G_z$ is a continuous fusion of homomorphisms defined on an open subgroup $H_0$ of $H$; Example 4.16  ibid.  shows that $\alpha$ can map nontrivially into each $G_z$. (This was perhaps surprising because one might reasonably, if naively,  have expected  from   the literature --- e.g., \cite{Coh, Il, Il-Sp, Il-Spr1, Il-St, Pham, Rud} --- that when a  homomorphism  $\vp = j_\alpha$ is completely positive,  $\alpha$  must map  homomorphically into a single subgroup of $\Delta(A)$.)  When $\vp = j_\alpha: A \ra B(H) $ is completely positive, the author of \cite{Sto21} was left with the question of whether $\alpha$ could map nontrivially into (a) infinitely many subgroups $G_z$ of $\Delta(A)$ and (b) into $G_z$, $G_w$ where $z, w$ are not comparable in $Z$. In the example below of a completely positive homomorphism $\vp = j_\alpha$, we also  observe that one can identify infinitely-many distinct infinite anti-chains $C$ in $(Z, \leq)$  such that  $\alpha$ maps nontrivially into each $G_z$ for $z$ in $C$. 
\er 

\bex \rm  \label{MainCPHomExample} Let $H$ be an infinite  group, $G = H^\N$ with the discrete topology. Consider the upper semilattice $(\F, \cup)$ where $\F = \{ A \subseteq \N: A \text{ is finite}\}$ and $A_\F^*(G) = A_{\B_\F}^*(G) =\ell^1-\bigoplus_{A \in \F} A_{\tau_A}(G)$, the GSA introduced in Example \ref{P(N) Example}. From Corollary \ref{Omnibus Example Cor} and Example  \ref{P(N) Example}, we have lower semilattice isomorphisms $$({\cal P}(\N), \cap) \cong (\HD(\F), \cap) \cong (\HD(\B_\F), \cap) \text{ given by }  B \mapsto \S_B = \{ A \in \F: A \subseteq B\} \mapsto \{ \tau_A : A \in \S_B\}$$ and 
$\Delta(A_\F^*(G)) = \bigcupdot_{\S \in \HD(\B_\F)} G_\S = \bigcupdot_{B \in {\cal P}(\N)} G_{\S_B},$ a semilattice of disjoint groups with compatible central identities, graded over the lower semilattice $({\cal P}(\N), \cap)$. Let $K$ be a nontrivial group and let $\alpha: K^\N \ra \Delta(A_\F^*(G))$ be the fusion of homomorphisms $\beta_B:K_B \ra G_{\S_B}$, compatible with $({\cal P}(\N), \cap) \cong (\HD(\B_\F), \cap )$, defined in Example  \ref{MainFusionExample}. Letting $K^\N$ have the discrete topology, $\vp = j_\alpha: A_\F^*(G) \ra B(K^\N)$ is a completely positive homomorphism by Theorem \ref{MainCPCCHomomThm}. As observed in Example \ref{MainFusionExample}, $$\alpha= \bigcupdot_{B\in{\cal P}(\N)} \beta_B \large{\vert}_{F_B} : K^\N = \bigcupdot_{B\in{\cal P}(\N)} F_B \ra \Delta(A_\F^*(G))  = \bigcupdot_{B\in{\cal P}(\N)} G_{\S_B} $$
where $F_B= \{ s\in K^\N: s_n = e_K \text{ for } n \in B, \ s_n \neq e_K \text{ for } n \in \N \bs B\}$, a nonempty subset of $K^\N$ for each $B$ in ${ \cal P}(\N)$. Thus, $\alpha$ maps nontrivially into $G_{\S_B}$ for each $B$ in ${\cal P}(\N)$. Clearly, $({\cal P}(\N), \cap)$ contains infinitely-many distinct infinite anti-chains and infinite chains, and, as noted in Example \ref{MainFusionExample}, $\S_B$ is nonprincipal in $\HD(\F)$ --- equivalently $\{\tau_A : A \in \S_B\}$ is nonprincipal in $\HD(\B_\F)$  --- for each of the uncountably-many infinite subsets $B$ of $\N$. 
\eex

A map $\alpha:E_0 \subseteq H \ra G_\B^* = \bigcupdot_{\S \in \HD(\B)} G_\S$ is  pw$^2$-affine (for piecewise-piecewise affine) if  $F_\S= F_\S^\alpha \in \Omega(H)$ and $\alpha_\S: = \alpha\vert_{F_\S} = \beta_S\vert_{F_\S}: F_\S \ra G_\S$ is piecewise affine for each $\S$ in $\HD(\B)$. This is consistent with the definition of pw$^2$-affine maps in \cite{Sto21}, (though therein such a map was required to map into only finitely many groups).  Since $\alpha = \bigcupdot_{\S \in \HD(\B)} \alpha_\S$, $j_\alpha u = \sum_{\S \in \HD(\B)} j_{\alpha_\S} u$ for $u \in A_\B^*(G)$.  When $G$ is amenable, Theorem \ref{CPCCCBHomomThm} with (\ref{KSbetaS General formula Eqn}) describe exactly when $\vp = j_\alpha: A_\B^*(G) \ra B(H)$ is a completely bounded homomorphism in terms of $\alpha$. We do not, in general, have a description of $\alpha$ for completely bounded homomorphisms $\vp = j_\alpha$ that is analogous to Theorem \ref{MainCPCCHomomThm}, but do have the following result.  

\bp  \rm Let $G$ be an amenable locally compact group, $\A=A_\B^*(G)$ a GSA. Consider the following conditions:  
\bi  \item[(a)] $\vp: A_\B^*(G) \ra B(H)$ is a completely bounded homomorphism; 
\item[(b)] there is an open subset $E_0$ of $H$ and a continuous pw$^2$-affine map $\alpha: E_0 \ra G_\B^* = \bigcupdot_{\S \in \HD(\B)} G_\S$ such that $\vp(u) = j_\alpha u $ for each $u$ in $A_\B^*(G)$.
\ei
If $\HD(\B) = \PHD(\B)$ and $U_\tau = \{ \sigma \in \B: \sigma \supsetneq \tau \}$ is finite for each $\tau \in \B$, then (a) implies (b). When $\B$ is finite, conditions (a) and (b) are equivalent and $E_0 \in \Omega(H)$. 
\ep 

\begin{proof}  Suppose that $\HD(\B) = \PHD(\B)$ and $U_\tau$ is finite for each $\tau \in \B$. Further, suppose that condition (a) holds and take $E_0$ and $\alpha: E_0 \ra G_\B^*$ satisfying condition (b) in Theorem \ref{CPCCCBHomomThm}.2. For $\S = \S_\tau$ in $\HD(\B)$, $F_{\S_\tau} = K_{\S_\tau} \bs \bigcup_{\sigma \in U_\tau} K_{\S_\tau}$ belongs to $\Omega(H)$, using Proposition \ref{Fusions of compatible maps Prop}.2. Since $\beta_\tau = \widetilde{\beta}_{\S_\tau}$ is piecewise affine and  $G_{\S_\tau} \ra G_\tau: s \mapsto \widetilde{s}$ is a topological group isomorphism,  $\beta_{\S_\tau}: K_{\S_\tau} \ra G_{\S_\tau}$   is also piecewise affine.  It is clear that the restriction of a piecewise affine map to a set in $\Omega(H)$ is piecewise affine, so each $\alpha_{\S_\tau} = \beta_{\S_\tau}\vert_{F_{\S_\tau}}$ is piecewise affine. Hence, $\alpha$ is pw$^2$-affine.  
If $\B$ is finite and $\S \in \HD(\B)$, then $\tau = \vee_\A \S \in \S$ and $\S = \S_\tau$, so $\HD(\B) = \PHD(\B)$. If, also,  (b) holds, then $\vp = j_\alpha: A_\B^*(G) \ra B(H)$ is a completely bounded homomorphism by \cite[Proposition 2.13]{Sto21}. (In this direction, amenability of $G$ is not needed.) 
\end{proof} 

\br  \label{Final Remark}  \rm 1.  Let  $({\cal C}, \cup)$ be the upper semilattice of cofinite subsets of $\N$ considered in Example \ref{CofiniteExample}. Then $A_{\B_{\cal C}}^*(G) = A_{\cal C}^*(G)$ is a GSA on a group $G$, that can be chosen to be abelian,  where $\B_{\cal C}$ is infinite but $\HD(\B_{\cal C}) = \PHD(\B_{\cal C})$ and $U_\tau$ is finite for each $\tau$ in $\B_{\cal C}$.  

\smallskip 

\noindent 2.  Some examples of groups $G$ for which $\TnqG$ is finite are found in Section 6 of \cite{Il-Spr1}.  By Theorem \ref{GSA MainThm2 Uniqueness of B}, for such groups $\B=\B_\A$ is finite for any GSA $\A$ on $G$. 
\er

\noindent {\sc Department of Pure Mathematics, University of
Waterloo, Waterloo ON, N2L 3G1, Canada;  email: {\tt nico.spronk@uwaterloo.ca}

\bigskip 

\noindent {\sc Department of Mathematics and Statistics, University
of Winnipeg, 515 Portage Avenue, Winnipeg, MB, Canada, R3B 2E9}; email: {\tt r.stokke@uwinnipeg.ca}

\bigskip 

\noindent {\sc Department of Mathematics, Pondicherry University, R.V.Nagar, Kalapet,
Puducherry – 605014}; email: {\tt aasaimanit@pondiuni.ac.in}


\begin{thebibliography}{99}


\bibitem{Ars} G. Arsac,  Sur l'espace de Banach engendr$\acute{\rm{e}}$
par les coefficients d'une repr$\acute{\rm{e}}$sentation unitaire,
{\it Publ. D$\acute{e}$p. Math. (Lyon)} 13 (1976), 1-101.

\bibitem{Ber-Jun-Mil} J.F. Berglund, H. Junghenn and P.  Milnes,  {\it Analysis on semigroups: function spaces, compactifications, representations,}  Canadian Mathematical Society Series of Monographs and Advanced Texts,  John Wiley \& Sons, Inc., New York, 1989.


\bibitem{Coh} P. J. Cohen, On homomorphisms of group algebras, {\it
Amer. J. Math} 82 (1960), 213-226.

\bibitem{Daws} M. Daws, Completely bounded homomorphisms of the Fourier algebra revisited, 
\it J. Group Theory \rm 25 (2022), no. 3, 579-600.

\bibitem{Eff-Rua} E.G. Effros and Z.-J. Ruan, {\it Operator Spaces}, Oxford
University Press, 2000.

\bibitem{Eym} P. Eymard,   L'alg\'ebre de Fourier d'un groupe localement compact,
 {\it Bull. Soc. Math. France}  92 (1964), 181-236.
 
 \bibitem{For} B. Forrest, {\it Fourier analysis on coset spaces},
Rocky Mountain J. Math. 28, (1998), no. 1, 173-189.

\bibitem{Gre}F.P. Greenleaf,  Norm decreasing homomorphisms of
 group algebras,  {\it Pacific J. Math.} 15 (1965), 1187-1219.
 
\bibitem{Il} M. Ilie,  On Fourier algebra homomorphisms, {\it J. Funct. Anal.}, 213
(2004), 88-110.


\bibitem{Il-Sp} M. Ilie and N. Spronk,  Completely bounded
homomorphisms of the Fourier algebra, {\it J. Funct. Anal.} 225
(2)(2005), 480-499.









\bibitem{Il-Spr1} M. Ilie and N. Spronk, The spine of a Fourier--Stieltjes algebra, {\it Proc. London Math. Soc.} (3) 94 (2007) 273-301.


\bibitem{Il-Spr2} M. Ilie and N. Spronk, Corrigendum: The spine of a Fourier--Stieltjes algebra, {\it Proc. London Math. Soc.} (3) 104 (2012) 859-863.

\bibitem{Il-St} M. Ilie and R. Stokke,  Weak$^*$-continuous homomorphisms of Fourier-Stieltjes algebras,
{\it  Math. Proc. Cambridge Philos. Soc.},  145 (2008),  107-120.

\bibitem{Inoue} J. Inoue, Some closed subalgebras of measure algebras and a generalization of P.J. Cohen's Theorem, {\it J. Math. Soc. Japan}, Vol. 23, No. 2 (1971), 278-294.

\bibitem{Kan} E. Kaniuth,  {\it A course in commutative Banach algebras}, Graduate Texts in Mathematics, Springer, New York, 2009.

\bibitem{Kan-Lau} E. Kaniuth and A.T.-M. Lau, \it Fourier and Fourier-Stieltjes algebras on locally compact groups, \rm  Mathematical Surveys and Monographs, 231, American Mathematical Society, Providence, RI, 2018.  

\bibitem{KSSY} M.E. Kroeker, A. Stephens,  R. Stokke,  R. Yee, Norm-multiplicative homomorphisms of
Beurling algebras, {\it J. Math. Anal. Appl.} 509 (2022), no. 1, Paper No. 125935, 25 pages.

\bibitem{Pau} V. Paulsen {\it Completely bounded maps and operator algebras},  Cambridge Studies in Advanced Mathematics, 78, Cambridge University Press, Cambridge, 2002.


\bibitem{Pham} H.L. Pham, Contractive homomorphisms of the
Fourier algebras, {\it Bull. London Math. Soc.} (2010) 42(5), 937-947.

\bibitem{Rud} W. Rudin, {\it Fourier analysis on groups}, Tracts in
Pure and Applied Mathematics, No. 12 Wiley  Interscience,  New
York-London 1962.




\bibitem{Spr}  N. Spronk, Weakly almost-periodic topologies, idempotents, and ideals, \it 
Indiana Univ. Math. J.  \rm 71 (2022), no. 6, 2671-2702.
\bibitem{Sto} R. Stokke, Spine-like subalgebras of Fourier--Stieltjes algebras, \it in preparation. \rm  


\bibitem{Sto11} R. Stokke, Homomorphisms of convolution algebras, {\it J. Funct. Anal.} 261 (2011), no. 12, 3665-3695 (2011).

\bibitem{Sto21} R. Stokke, Homomorphisms of Fourier--Stieltjes algebras, \it
Studia Math. \rm 258 (2021), no. 2, 175-220.

\bibitem{Tak} M. Takesaki, \it Theory of Operator Algebras I, \rm Encyclopedia of Mathematical Sciences Vol. 124, Springer--Verlag Berlin Heidelberg, 2002. 

\bibitem{Tay} J. L. Taylor, \it Measure algebras, \rm CBMS Regional Conference Series in Mathematics 16 (American Mathematical Society, Providence, RI, 1973).

 \bibitem{Tha} A. Thamizhazhagan, On the structure of invertible elements in certain Fourier--Stieltjes algebras, \it Studia Math. \rm 257 (2021), no. 3, 347-360.

\bibitem{Wal1} M. Walter,  $W\sp{*} $-algebras and nonabelian harmonic
analysis, {\it  J. Funct. Anal.} 11 (1972), 17--38.

\bibitem{Wal2} M. Walter, 
On the structure of the Fourier-Stieltjes algebra, {\it
Pacific J. Math.} 58 (1975), no. 1, 267-281. 


\end{thebibliography}
\end{document}